\documentstyle[12pt]{article}
\def\Bbb R{{\rm \bf R}}
\def\proclaim#1{\vskip2mm{\bf #1}\em}
\def\endproclaim{\em \vskip2mm}
\def\tag#1{\eqno(#1)}
\def\gathered{\begin{array}{c}}
\def\endgathered{\end{array}}
\def\text{\mbox}

\begin{document}

\title {Integrating the probe and singular sources methods:III.  Mixed obstacle case}
\author{Masaru IKEHATA\footnote{
Laboratory of Mathematics,
Graduate School of Advanced Science and Engineering,
Hiroshima University,
Higashihiroshima 739-8527, JAPAN.
e-mail address: ikehataprobe@gmail.com}
\footnote{Professor Emeritus at Gunma University}
\footnote{Professor Emeritus at Hiroshima University}
}
\maketitle

\begin{abstract}
The main purpose of this paper is to develop further the integrated theory of the probe and singular sources methods (IPS) which may work
for a group of inverse obstacle problems.
Here as a representative and typical member of the group, an inverse obstacle problem governed by the Helmholtz equation with 
a fixed wave number in a bounded domain is considered.  It is assumed that the solutions of the Helmholtz equation
outside the set of unknown obstacles satisfy the homogeneous Dirichlet or Neumann boundary conditions
on each surface of obstacles.   This is the case when two extreme types of obstacles are embedded in a medium.
By considering this case, not only a concise technique for IPS is introduced but also
a general correspondence principle from IPS to the probe method is suggested.
Besides, as a corollary it is shown that the probe method together with
the  singular sources method reformulated in terms of the probe method has
the Side B under a smallness conditions on the wave number $k$, which is the blowing up property of a sequence
computed from the associated Dirichlet-to-Neumann map.

\noindent
AMS:  35R30, 78A46, 35J05, 35E05

\end{abstract}



\section{Introduction}

Both the probe method of Ikehata \cite{IProbe,I2} (later reformulated in \cite{INew})
and singular sources method of Potthast 
\cite{P1,NP} now become well-known classical
analytical methods for reconstruction issue of {\it inverse obstacle problems} governed by partial differential equations.
This paper is concerned with the {\it integrated theory} of the probe and singular sources methods (IPS),
which is initiated by the author himself  in \cite{IPS,IPS2}.
In particular, we focus on the role of IPS in deriving the probe and singular sources methods
together with introducing  a technique to treat some kind of inverse obstacle problems 
governed by partial differential equations.  For the purpose we consider a prototype inverse obstacle problem
governed by the Helmholtz equation with a fixed wave number.

Now let us formulate the prototype problem.
Let $\Omega$ be a bounded domain of $\Bbb R^3$ with Lipschitz-boundary \cite{Gr}.
We denote by $D$ a mathematical model of {\it discontinuity}
embedded in the background medium
$\Omega$.  
We assume that $D$ takes the form $D=D_n\cup D_d$, where $D_n$ and $D_d$ are open subsets of $\Bbb R^3$ with Lipschitz-boundary
with $\overline{D_n}\cap\overline{D_d}=\emptyset$, $\overline{D_n}
\cup\overline{D_d}\subset\Omega$ and that $\Omega\setminus(
\overline{D_n}\cup\overline{D_d})$ is connected.  We denote by $\nu$ the unit outward normal vector to not only $\partial\Omega$ but also
$\partial D$.
On the surfaces of $D_n$ and $D_d$ two boundary conditions of different type are imposed
as specified below.

Let $k\ge 0$.  Given an arbitrary $f\in H^{\frac{1}{2}}(\partial\Omega)$, let $u=u(x)$ in $H^1(\Omega\setminus\overline{D})$ be the weak solution of
$$
\left\{
\begin{array}{ll}
\displaystyle
\Delta u+k^2u=0, & x\in\Omega\setminus\overline{D},
\\
\\
\displaystyle
\frac{\partial u}{\partial\nu}=0, & x\in\partial D_n,
\\
\\
\displaystyle
u=0, & x\in\partial D_d,
\\
\\
\displaystyle
u=f, & x\in\partial\Omega.
\end{array}
\right.
\tag {1.1}
$$
This means that, $u=f$ on $\partial\Omega$, $u=0$ on $\partial D_d$ in the sense of the trace
and, for all $\varphi\in H^1(\Omega\setminus\overline{D})$ with $\varphi=0$ on $\partial\Omega$ and $\varphi=0$ on $\partial D_d$ in the sense of the trace, we have
$$\displaystyle
-\int_{\Omega\setminus\overline{D}}\,\nabla u\cdot\nabla\varphi\,dx+\int_{\Omega\setminus\overline{D}}k^2u\varphi\,dx=0.
$$
Then the bounded linear functional $\frac{\partial u}{\partial\nu}\vert_{\partial\Omega}\in H^{-\frac{1}{2}}(\partial\Omega)$
is well defined via the formula
$$\begin{array}{ll}
\displaystyle
\left<\frac{\partial u}{\partial\nu}\vert_{\partial\Omega}, g\right>=\int_{\Omega\setminus\overline{D}}\,\nabla u\cdot\nabla\phi\, dx-\int_{\Omega\setminus\overline{D}}k^2u\phi\,dx,
&
\displaystyle
g\in H^{\frac{1}{2}}(\partial\Omega),
\end{array}
$$
where $\phi\in H^1(\Omega\setminus\overline{D})$ such that $\phi=g$ on $\partial\Omega$ and $\phi=0$ on $\partial D_d$ in the sense of the trace.
Note that unless otherwise speciefied the functions appearing in this paper are always real-valued;
the symbol $\nu$ denotes the unit outwrd normal vector field on $\partial\Omega$ and $\partial D=\partial D_n\cup\partial D_d$.

In this paper, by considering the prototype inverse obstacle problem  mentioned below, we further
develop a technique to the integrated theory of the probe and singular sources methods.

$\quad$

{\noindent\bf Problem.}  Extract information about the geometry of $D_n$ and $D_d$ from the 
$\frac{\partial u}{\partial\nu}\vert_{\partial\Omega}$ corresponding to {\it infinitely many} $f$.

$\quad$

For the problem to have meaning we impose a restriction on $k$:

$\quad$

{\bf\noindent Assumption 1.}  The boundary value problem (1.1)  with $f=0$ has only a trivial solution.

$\quad$

\noindent
Under Assumption 1  it is well known that the weak solution $u$ of  (1.1) exists and unique. 
Then the map
$$\displaystyle
\Lambda_D:H^{\frac{1}{2}}(\partial\Omega)\rightarrow H^{-\frac{1}{2}}(\partial\Omega),
$$
is well-defined by 
$$
\Lambda_Df=\frac{\partial u}{\partial\nu}\vert_{\partial\Omega}.
$$
This is called the Dirichlet-to-Neumann map.  
So Problem becomes the extraction problem of information about the geometry 
of $D_n$ and $D_d$ from the graph of the Dirichlet-to-Neumann map $\Lambda_D$ or its  partial knowledge.

It follows from the definition we have the symmetry:
for all $f\in H^{\frac{1}{2}}(\partial\Omega)$ and $g\in H^{\frac{1}{2}}(\partial\Omega)$
$$\displaystyle
<\Lambda_D f,g>=<\Lambda_Dg,f>.
$$
Besides, if $f\in H^{\frac{3}{2}}(\partial\Omega)$
and both $\partial\Omega$ and $\partial D$ are $C^{2}$, then $u\in H^2(\Omega\setminus\overline{D})$
and thus $\Lambda_Df=\frac{\partial u}{\partial\nu}\vert_{\partial\Omega}\in H^{\frac{1}{2}}(\partial\Omega)$
in the sense of the trace \cite{Gr}.  Then,  integration by parts (e.g, Lemma 1.5.3.7 of \cite{Gr}) yields the surface integral expression of $<\Lambda_Df,g>$ for
all  $f\in H^{\frac{3}{2}}(\partial\Omega)$ and $g\in H^{\frac{1}{2}}(\partial\Omega)$
$$\displaystyle
<\Lambda_D f,g>
=\int_{\partial\Omega}\,\Lambda_D f(z)g(z)\,dS(z).
$$
In this paper, we always consider $k$ such that Assumption 1 is satisfied and, unless otherwise stated the $C^2$-regularity
of $\partial\Omega$ and $\partial D$ are assumed.  Note that $k=0$ satisfies Assumption 1

In \cite{IPS} by considering the case that $D_d=\emptyset$ in (1.1) and $k=0$, the author introduced the integrated theory
of the probe and singular sources methods.  In \cite{IPS2} IPS has been applied to
an inverse obstacle problem governed by the Stokes system.  Therein a technique to treat a system is introduced.
In this paper we pursuit IPS further by considering the case when $D_d\not=\emptyset$ and $k\ge 0$.  
Especially, with application to elastic bodies in mind, it would be interesting to consider such a case.
It is expected that this situation causes some problems due to the coexistence of two different boundary conditions
and $k\not=0$.
Besides, it cannot be said that the IPS concept was thoroughly developed in \cite{IPS,IPS2} in the sense that
the treatment of the probe method therein is independent from IPS.
This time we would like to show that not only the singular sources method
but also the probe method itself is derived from IPS.

It should be pointed out that, using the {\it original} probe method \cite{Iwave}, this type of problem itself has been considered by Cheng-Liu-Nakamura-Wang \cite{CLNW}.  However, they do not have the view point of the IPS developed in this paper.  See also Remark 3.4 for more detailed comparison.
Note that, for mixed obstacles placed in the {\it whole space} there are some applications of 
the {\it factorization method} \cite{KG}, \cite{L}, \cite{F}, {\it monotonicity method} \cite{AG} and both \cite{F2}.

\subsection{The IPS function}

The IPS in this paper starts with introducing 
a family of singular solutions for the back ground medium.

Let ${\cal G}=\{G(\,\cdot\,,x)\}_{x\in\Omega}$ be a family of distributions in $\Omega$ indexed
with $x\in\Omega$ having the form
$$\displaystyle
G(y,x)=G(y-x)+H(y,x),
\tag {1.2}
$$
where 
$$\displaystyle
G(y-x)=\frac{\cos k\vert y-x\vert}{4\pi\vert y-x\vert}
$$
and $H(\,\cdot\,,x)\in H^2(\Omega)$ is a real-valued solution of the Helmholtz equation in $\Omega$ 
such that, for each $\epsilon>0$
$$\displaystyle
\sup_{x\in\Omega,\,\text{dist}\,(x,\partial\Omega)>\epsilon\,}\Vert H(\,\cdot\,,x)\Vert_{H^2(\Omega)}<\infty.
\tag {1.3}
$$
Note that $G(\,\cdot\,-x)$ coincides with the real part of the standard (complex-valued) fundamental solution
of the Helmholtz equation 
$$\displaystyle
\Phi(y-x)=\frac{e^{ik\vert y-x\vert}}{4\pi\vert y-x\vert}.
$$
Since the imaginary part of  $\Phi(\,\cdot\,-x)$ has the unique extension to the whole space as the entire 
solution of the Helmholtz equation, the function $G(y-x)$ also satisfies 
$$\displaystyle
\Delta G(\,\cdot\,-x)+k^2 G(\,\cdot\,-x)+\delta(\,\cdot\,-x)=0
$$
as the distribution of $y\in\Bbb R^3$ for each fixed $x\in\Bbb R^3$.

$\quad$

{\bf\noindent Definition 1.1.}
Given ${\cal G}$ and $x\in\Omega\setminus\overline{D}$ let $W=W_x(y;{\cal G})=W(y)$ in $H^2(\Omega\setminus\overline{D})$ be the solution of 
$$\left\{
\begin{array}{ll}
\displaystyle
\Delta W+k^2W=0, & y\in\Omega\setminus\overline{D},\\
\\
\displaystyle
\frac{\partial W}{\partial\nu}=-\frac{\partial}{\partial\nu}G(y,x), & y\in\partial D_n,\\
\\
\displaystyle
W=-G(y,x), & y\in\partial D_d,\\
\\
\displaystyle
W=G(y,x), & y\in\partial\Omega.
\end{array}
\right.
\tag {1.4}
$$
We call the function $\Omega\setminus\overline{D}\ni x\longmapsto W_x(x;{\cal G})$ the IPS function
based on ${\cal G}$ for obstacle $D$.

$\quad$

\noindent
Hereafter we simply write $W_x(y;{\cal G})=W_x(y)$.
Needless to say, Assumption 1 ensures the uniqe solvability of  equations (1.4) in the class $H^2(\Omega\setminus\overline{D})$.
The sudden appearance of the system (1.4) seems strange, however, 
it is a natural extension of  the corresponding one firstly introduced in \cite{IPS} in the case when $D_d=\emptyset$
and $k=0$.
Besides, we will see in Section 3 that the IPS function 
generates the indicator function (see Definition 3.2) for the probe method.

Since the system (1.4) is linear, we have the natural and trivial decomposition
of the solution as 
$$\begin{array}{ll}
\displaystyle
W_x(y)=w_x(y)+w_x^1(y),
&
y\in\Omega\setminus\overline{D},
\end{array}
\tag {1.5}
$$
where $w_x=w_x(y;{\cal G})=w(y)$ in $H^2(\Omega\setminus\overline{D})$ solves
$$\left\{
\begin{array}{ll}
\displaystyle
\Delta w+k^2w=0, & y\in\Omega\setminus\overline{D},\\
\\
\displaystyle
\frac{\partial w}{\partial\nu}=-\frac{\partial}{\partial\nu}G(y,x), & y\in\partial D_n,\\
\\
\displaystyle
w=-G(y,x), & y\in\partial D_d,\\
\\
\displaystyle
w=0, & y\in\partial\Omega
\end{array}
\right.
\tag {1.6}
$$
and $w_x^1=w_x^1(y;{\cal G})=w(y)$ in $H^2(\Omega\setminus\overline{D})$ solves
$$\left\{
\begin{array}{ll}
\displaystyle
\Delta w+k^2w=0, & y\in\Omega\setminus\overline{D},\\
\\
\displaystyle
\frac{\partial w}{\partial\nu}=0, & y\in\partial D_n,\\
\\
\displaystyle
w=0, & y\in\partial D_d,\\
\\
\displaystyle
w=G(y,x), & y\in\partial\Omega.
\end{array}
\right.
\tag {1.7}
$$
The unique solvability of the boundary value problems (1.6) and (1.7) are also a consequence of Assumption 1.

Needless to say, from (1.5) we have
$$\begin{array}{ll}
\displaystyle
W_x(x)=w_x(x)+w_x^1(x), & x\in\Omega\setminus\overline{D}.
\end{array}
\tag {1.8}
$$
We call this the {\it outer decomposition} or natural decomposition of IPS function.

The IPS for Problem is based on the discovery of the following two representation formulae.

\proclaim{\noindent Theorem 1.1.} 
Let $x\in\Omega\setminus\overline{D}$.

\noindent
(i)  We have the expression focused on the Neumann obstacle
$$\begin{array}{ll}
\displaystyle
W_x(x)
&
\displaystyle
=\Vert\nabla G(\,\cdot\,,x)\Vert_{L^2(D_n)}^2-k^2\Vert G(\,\cdot\,,x)\Vert_{L^2(D_n)}^2
\\
\\
\displaystyle
&
\displaystyle
\,\,\,
+\Vert \nabla(w_x+(\epsilon_x)_n)\Vert_{L^2(\Omega\setminus\overline{D})}^2
-k^2\Vert w_x+(\epsilon_x)_n\Vert_{L^2(\Omega\setminus\overline{D})}^2
\\
\\
\displaystyle
&
\displaystyle
\,\,\,
-\int_{\partial\Omega}\,G(z,x)\frac{\partial}{\partial\nu}G(z,x)\,dS(z)
\\
\\
\displaystyle
&
\displaystyle
\,\,\,
+\Vert \nabla w_x^1\Vert_{L^2(\Omega\setminus\overline{D})}^2-k^2\Vert w_x^1\Vert_{L^2(\Omega\setminus\overline{D})}^2
\\
\\
\displaystyle
&
\displaystyle
\,\,\,
-\Vert \nabla(\epsilon_x)_n\Vert_{L^2(\Omega\setminus\overline{D})}^2
+k^2\Vert(\epsilon_x)_n\Vert_{L^2(\Omega\setminus\overline{D})}^2
\\
\\
\displaystyle
&
\displaystyle
\,\,\,
-\Vert\nabla G(\,\cdot\,,x)\Vert_{L^2(D_d)}^2+k^2\Vert G(\,\cdot\,,x)\Vert_{L^2(D_d)}^2,
\end{array}
\tag {1.9}
$$
where $(\epsilon_x)_n=(\epsilon_x)_n(y;{\cal G})=\epsilon(y)$ in $H^2(\Omega\setminus\overline{D})$ solves
$$\left\{
\begin{array}{ll}
\displaystyle
\Delta \epsilon+k^2\epsilon=0, & y\in\Omega\setminus\overline{D},\\
\\
\displaystyle
\epsilon=0, & y\in\partial D_n,\\
\\
\displaystyle
\epsilon=G(y,x), & y\in\partial D_d,\\
\\
\displaystyle
\epsilon=0, & y\in\partial\Omega.
\end{array}
\right.
\tag {1.10}
$$
(ii)  We have the expression focused on the Dirichlet obstacle
$$\begin{array}{ll}
\displaystyle
W_x(x)
&
\displaystyle
=-\Vert\nabla G(\,\cdot\,,x)\Vert_{L^2(D_d)}^2+k^2\Vert G(\,\cdot\,,x)\Vert_{L^2(D_d)}^2
\\
\\
\displaystyle
&
\displaystyle
\,\,\,
-\Vert \nabla(w_x+(\epsilon_x)_d)\Vert_{L^2(\Omega\setminus\overline{D})}^2
+k^2\Vert w_x+(\epsilon_x)_d\Vert_{L^2(\Omega\setminus\overline{D})}^2
\\
\\
\displaystyle
&
\displaystyle
\,\,\,
-\int_{\partial\Omega}\,G(z,x)\frac{\partial}{\partial\nu}G(z,x)\,dS(z)
\\
\\
\displaystyle
&
\displaystyle
\,\,\,
+\Vert \nabla w_x^1\Vert_{L^2(\Omega\setminus\overline{D})}^2-k^2\Vert w_x^1\Vert_{L^2(\Omega\setminus\overline{D})}^2
\\
\\
\displaystyle
&
\displaystyle
\,\,\,
+\Vert \nabla(\epsilon_x)_d\Vert_{L^2(\Omega\setminus\overline{D})}^2-k^2\Vert (\epsilon_x)_d\Vert_{L^2(\Omega\setminus\overline{D})}^2
\\
\\
\displaystyle
&
\displaystyle
\,\,\,
+\Vert\nabla G(\,\cdot\,,x)\Vert_{L^2(D_n)}^2-k^2\Vert G(\,\cdot\,,x)\Vert_{L^2(D_n)}^2,
\end{array}
\tag {1.11}
$$
where $(\epsilon_x)_d=(\epsilon_x)_d(y;{\cal G})=\epsilon(y)$ in $H^2(\Omega\setminus\overline{D})$ solves
$$\left\{
\begin{array}{ll}
\displaystyle
\Delta\epsilon+k^2\epsilon=0, & y\in\Omega\setminus\overline{D},\\
\\
\displaystyle
\frac{\partial\epsilon}{\partial\nu}=\frac{\partial}{\partial\nu}G(y,x), & y\in\partial D_n,\\
\\
\displaystyle
\frac{\partial\epsilon}{\partial\nu}=0, & y\in\partial D_d,\\
\\
\displaystyle
\epsilon=0, & y\in\partial\Omega.
\end{array}
\right.
\tag {1.12}
$$

\endproclaim

\noindent
Introducing the functions $(\epsilon_x)_n$ and $(\epsilon_x)_d$ helps us to
write the IPS function in terms of energy integrals (1.9) and (1.11).  
We call this technique the {\it method of complementing function}.
A clear advantage is that: roughly speaking, we can immediately see that $(\epsilon_x)_n$ is bounded in
$H^2(\Omega\setminus\overline{D})$ if $x$ is close to a point on $\partial D_n$;
 $(\epsilon_x)_d$ is bounded in
$H^2(\Omega\setminus\overline{D})$ if $x$ is close to a point on $\partial D_d$.

It shoulde be pointed out that the two expressions (1.9) and (1.11) contain the common terms of  integrals
$$\begin{array}{l}
\displaystyle
\,\,\,\,\,\,
-\int_{\partial\Omega}\,G(z,x)\frac{\partial}{\partial\nu}G(z,x)\,dS(z)
+\Vert \nabla w_x^1\Vert_{L^2(\Omega\setminus\overline{D})}^2-k^2\Vert w_x^1\Vert_{L^2(\Omega\setminus\overline{D})}^2.
\end{array}
$$
Except for those, the expression of the right-hand side on (1.11) coincides with the one on (1.9) mutiplied by $(-1)$ and replaced $(n,d)$ with $(d,n)$.

The following corollary is a direct consequence of the facts listed below:

\noindent
(a)  the well-posedness of the boundary value problems (1.6), (1.7), (1.10) and (1.12);

\noindent
(b)   the expressions (1.9) and (1.11);

\noindent
 (c)  the property that for any finite cone $V$ with vertex at the origin
 $x=0$
$$\displaystyle
\int_V\vert\nabla G(z-x)\vert^2\,dz=\infty.
$$

\proclaim{\noindent Corollary 1.1.}
The IPS function $W_x(x)$ satisfies (i),  (ii) and (iii) listed below.

\noindent
(i)  $\lim_{x\rightarrow a\in\partial D_n}W_x(x)=\infty$.

\noindent
(ii)  $\lim_{x\rightarrow b\in\partial D_d}W_x(x)=-\infty$.

\noindent
(iii)  For each $\epsilon_i>0$, $i=1,2$ 
$$\displaystyle
\sup_{x\in\Omega\setminus\overline{D},\,\text{dist}\,(x,\partial D)>\epsilon_1,\,\text{dist}\,(x,\partial\Omega)>\epsilon_2}\,\vert W_x(x)\vert<\infty.
$$

\endproclaim

{\it\noindent Proof.}
In what follows we denote by $C_1, C_2,\cdots$ positive numers independent of $x$.
Using (1.9), from (a) and (b) together with (1.2) and (1.3) we have:
as $x\rightarrow a\in\partial D_n$
$$
\displaystyle
W_x(x)
\ge \Vert\nabla G(\,\cdot\,,x)\Vert_{L^2(D_n)}^2-2k^2\Vert w\Vert_{L^2(\Omega\setminus\overline{D})}^2-C_1.
$$
Here, by Lemma B in Appendix we have, as $x\rightarrow a\in\partial D_n$,
$$
\displaystyle
W_x(x)
\ge \Vert\nabla G(\,\cdot\,,x)\Vert_{L^2(D_n)}^2-C_1-C_2.
$$
Then (c) together with Fatou's lemma yields (i).  
Next consider the case when $x\rightarrow b\in\partial D_d$.
It follows from (1.11)
$$
\displaystyle
-W_x(x)
\ge \Vert\nabla G(\,\cdot\,,x)\Vert_{L^2(D_d)}^2-2k^2\Vert w\Vert_{L^2(\Omega\setminus\overline{D})}^2-C_3.
$$
Again, Lemma B yields, as $x\rightarrow b\in\partial D_d$
$$\displaystyle
-W_x(x)
\ge C_4\Vert\nabla G(\,\cdot\,,x)\Vert_{L^2(D_d)}^2-C_3-C_5.
$$
This yields the validity of (ii).
The validity of statement (iii) is almost clear.

\noindent
$\Box$

\noindent
Therefore, using the asymptotic behaviour of  IPS function $W_x(x)$ as $x$ approches $\partial D$, one can distinguish the soft obstacle $D_d$ and hard obstacle $D_n$.  In particular, we know that IPS function does not have a definite sign unlike a single type of obstacle case \cite{IPS}.

$\quad$

{\bf\noindent Remark 1.1.}  In this paper, about the choice of  the family $\{H(\,\cdot\,,x)\}_{x\in\Omega}$ in (1.2)
we consider only the two cases.
The first is the case when $H(y,x)\equiv 0$.   
In this case we denote ${\cal G}$ by ${\cal G}^0$.
Then $G(\,\cdot\,,x)$ coincides with $G(\,\cdot\,-x)$.
The second is:  we impose the boundary condition
$$\begin{array}{ll}
\displaystyle
H(y,x)=-\frac{\cos k\vert y-x\vert}{4\pi\vert y-x\vert}, & y\in\partial\Omega.
\end{array}
\tag {1.13}
$$
Under the assumption that $k^2$ is not a Dirichlet eigenvalue for the minus Laplacian $-\Delta$ in $\Omega$, for each $x\in\Omega$
the $H(y,x)$ exists and is unique,  and satisfies (1.3).
The function $G(\,\cdot\,,x)$ is nothing but the Green function for the domain $\Omega$ with the source point at $x\in\Omega$. 
In this case we denote ${\cal G}$ by ${\cal G}^*$.   Then $W_x=0$ on $\partial\Omega$ for $x\in\Omega\setminus\overline{D}$ and $w_x^1\equiv 0$.
Hereafter unless otherwise stated, we always impose the condition on $k^2$ mentioned above when considering ${\cal G}^*$.

Theorem 1.1 yields the following corollary.

\proclaim{\noindent Corollary 1.2.}   Choose ${\cal G}={\cal G}^0$.  Let $w_x=w_x(\,\cdot\,;{\cal G}^0)$ and $w_x^1=w_x^1(\,\cdot\,;{\cal G}^0)$.

\noindent
(i)  We have the expression focused on the Neumann obstacle
$$\begin{array}{ll}
\displaystyle
W_x(x)
&
\displaystyle
=\Vert\nabla G(\,\cdot\,-x)\Vert_{L^2(D_n)}^2
-k^2\Vert G(\,\cdot\,-x)\Vert_{L^2(D_n)}^2
\\
\\
\displaystyle
&
\displaystyle
\,\,\,
+\Vert \nabla(w_x+(\epsilon_x)_n)\Vert_{L^2(\Omega\setminus\overline{D})}^2-k^2\Vert w_x+(\epsilon_x)_n\Vert_{L^2(\Omega\setminus\overline{D})}^2
\\
\\
\displaystyle
&
\displaystyle
\,\,\,
-\int_{\partial\Omega}\,G(z-x)\frac{\partial}{\partial\nu}G(z-x)\,dS(z)
\\
\\
\displaystyle
&
\displaystyle
\,\,\,
+\Vert \nabla w_x^1\Vert_{L^2(\Omega\setminus\overline{D})}^2-k^2\Vert w_x^1\Vert_{L^2(\Omega\setminus\overline{D})}^2
\\
\\
\displaystyle
&
\displaystyle
\,\,\,
-\Vert \nabla(\epsilon_x)_n\Vert_{L^2(\Omega\setminus\overline{D})}^2
+k^2\Vert (\epsilon_x)_n\Vert_{L^2(\Omega\setminus\overline{D})}^2
\\
\\
\displaystyle
&
\displaystyle
\,\,\,
-\Vert\nabla G(\,\cdot\,-x)\Vert_{L^2(D_d)}^2+k^2\Vert G(\,\cdot\,-x)\Vert_{L^2(D_d)}^2,
\end{array}
$$
where $(\epsilon_x)_n=(\epsilon_x)_n(\,\cdot\,;{\cal G}^0)$.

\noindent
(ii)  We have the expression focused on the Dirichlet obstacle
$$\begin{array}{ll}
\displaystyle
W_x(x)
&
\displaystyle
=-\Vert\nabla G(\,\cdot\,-x)\Vert_{L^2(D_d)}^2+k^2\Vert G(\,\cdot\,-x)\Vert_{L^2(D_d)}^2
\\
\\
\displaystyle
&
\displaystyle
\,\,\,
-\Vert \nabla(w_x+(\epsilon_x)_d)\Vert_{L^2(\Omega\setminus\overline{D})}^2+k^2\Vert (w_x+(\epsilon_x)_d)\Vert_{L^2(\Omega\setminus\overline{D})}^2
\\
\\
\displaystyle
&
\displaystyle
\,\,\,
-\int_{\partial\Omega}\,G(z-x)\frac{\partial}{\partial\nu}G(z-x)\,dS(z)
\\
\\
\displaystyle
&
\displaystyle
\,\,\,
+\Vert \nabla w_x^1\Vert_{L^2(\Omega\setminus\overline{D})}^2-k^2\Vert w_x^1\Vert_{L^2(\Omega\setminus\overline{D})}^2
\\
\\
\displaystyle
&
\displaystyle
\,\,\,
+\Vert \nabla(\epsilon_x)_d\Vert_{L^2(\Omega\setminus\overline{D})}^2
-k^2\Vert (\epsilon_x)_d\Vert_{L^2(\Omega\setminus\overline{D})}^2
\\
\\
\displaystyle
&
\displaystyle
\,\,\,
+\Vert\nabla G(\,\cdot\,-x)\Vert_{L^2(D_n)}^2-k^2\Vert G(\,\cdot\,-x)\Vert_{L^2(D_n)}^2,
\end{array}
$$
where $(\epsilon_x)_d=(\epsilon_x)_d(\,\cdot\,;{\cal G}^0)$.

\endproclaim

$\quad$

{\bf\noindent Remark 1.1.}
In particular, if $k=0$, then for all $x\in\Omega$ and $y\in\Omega$ one can rewrite
$$\displaystyle
-\int_{\partial\Omega}\,G(z-x)\frac{\partial}{\partial\nu}G(z-y)\,dS(z)=
\int_{\Bbb R^3\setminus\overline{\Omega}}\,\nabla G(z-x)\cdot\nabla G(z-y)\,dz.
\tag {1.14}
$$
Note that the integrand of this right-hand side is absolutely integrable.
Thus, the formulae in Corollary 1.2 become
$$\begin{array}{ll}
\displaystyle
W_x(x)
&
\displaystyle
=\Vert\nabla G(\,\cdot\,-x)\Vert_{L^2(D_n)}^2
+\Vert \nabla(w_x+(\epsilon_x)_n)\Vert_{L^2(\Omega\setminus\overline{D})}^2
\\
\\
\displaystyle
&
\displaystyle
\,\,\,
+\Vert\nabla G(\,\cdot\,-x)\Vert_{L^2(\Bbb R^3\setminus\overline{\Omega})}^2
\\
\\
\displaystyle
&
\displaystyle
\,\,\,
+\Vert \nabla w_x^1\Vert_{L^2(\Omega\setminus\overline{D})}^2
-\Vert \nabla(\epsilon_x)_n\Vert_{L^2(\Omega\setminus\overline{D})}^2
\\
\\
\displaystyle
&
\displaystyle
\,\,\,
-\Vert\nabla G(\,\cdot\,-x)\Vert_{L^2(D_d)}^2
\end{array}
$$
and
$$\begin{array}{ll}
\displaystyle
W_x(x)
&
\displaystyle
=-\Vert\nabla G(\,\cdot\,-x)\Vert_{L^2(D_d)}^2
-\Vert \nabla(w_x+(\epsilon_x)_d)\Vert_{L^2(\Omega\setminus\overline{D})}^2
\\
\\
\displaystyle
&
\displaystyle
\,\,\,
+\Vert\nabla G(\,\cdot\,-x)\Vert_{L^2(\Bbb R^3\setminus\overline{\Omega})}^2
\\
\\
\displaystyle
&
\displaystyle
\,\,\,
+\Vert \nabla w_x^1\Vert_{L^2(\Omega\setminus\overline{D})}^2
+\Vert \nabla(\epsilon_x)_d\Vert_{L^2(\Omega\setminus\overline{D})}^2
\\
\\
\displaystyle
&
\displaystyle
\,\,\,
+\Vert\nabla G(\,\cdot\,-x)\Vert_{L^2(D_n)}^2.
\end{array}
$$

$\quad$

\noindent
And also we have

\proclaim{\noindent Corollary 1.3.}   Choose ${\cal G}={\cal G}^*$.  Then $w_x^1(\,\cdot\,;{\cal G}^*)=0$ and $W_x(\,\cdot\,;{\cal G}^*)=w_x(\,\cdot\,;{\cal G}^*)$.

\noindent
(i)  We have the expression focused on the Neumann obstacle
$$\begin{array}{ll}
\displaystyle
W_x(x)
&
\displaystyle
=\Vert\nabla G(\,\cdot\,,x)\Vert_{L^2(D_n)}^2-k^2\Vert G(\,\cdot\,,x)\Vert_{L^2(D_n)}^2
\\
\\
\displaystyle
&
\displaystyle
\,\,\,
+\Vert \nabla(w_x+(\epsilon_x)_n)\Vert_{L^2(\Omega\setminus\overline{D})}^2
-k^2\Vert (w_x+(\epsilon_x)_n)\Vert_{L^2(\Omega\setminus\overline{D})}^2
\\
\\
\displaystyle
&
\displaystyle
\,\,\,
-\Vert \nabla(\epsilon_x)_n\Vert_{L^2(\Omega\setminus\overline{D})}^2+k^2\Vert (\epsilon_x)_n\Vert_{L^2(\Omega\setminus\overline{D})}^2
\\
\\
\displaystyle
&
\displaystyle
\,\,\,
-\Vert\nabla G(\,\cdot\,,x)\Vert_{L^2(D_d)}^2+k^2\Vert G(\,\cdot\,,x)\Vert_{L^2(D_d)}^2,
\end{array}
$$
where $(\epsilon_x)_n=(\epsilon_x)_n(\,\cdot\,;{\cal G}^*)$.

\noindent
(ii)  We have the expression focused on the Dirichlet obstacle
$$\begin{array}{ll}
\displaystyle
W_x(x)
&
\displaystyle
=-\Vert\nabla G(\,\cdot\,,x)\Vert_{L^2(D_d)}^2+k^2\Vert G(\,\cdot\,,x)\Vert_{L^2(D_d)}^2
\\
\\
\displaystyle
&
\displaystyle
\,\,\,
-\Vert \nabla(w_x+(\epsilon_x)_d)\Vert_{L^2(\Omega\setminus\overline{D})}^2+k^2\Vert (w_x+(\epsilon_x)_d)\Vert_{L^2(\Omega\setminus\overline{D})}^2
\\
\\
\displaystyle
&
\displaystyle
\,\,\,
+\Vert \nabla(\epsilon_x)_d\Vert_{L^2(\Omega\setminus\overline{D})}^2-k^2\Vert (\epsilon_x)_d\Vert_{L^2(\Omega\setminus\overline{D})}^2
\\
\\
\displaystyle
&
\displaystyle
\,\,\,
+\Vert\nabla G(\,\cdot\,,x)\Vert_{L^2(D_n)}^2-k^2\Vert G(\,\cdot\,,x)\Vert_{L^2(D_n)}^2,
\end{array}
$$
where $(\epsilon_x)_d=(\epsilon_x)_d(\,\cdot\,;{\cal G}^*)$.

\endproclaim

$\quad$

\subsection{IPS function to DN map}

In this section, we assume that Assumption 1 for the case $D=\emptyset$ is also satisfied.
We denote $\Lambda_D$ by $\Lambda_0$ if $D=\emptyset$.
In the probe method the form
$$\displaystyle
<(\Lambda_0-\Lambda_D)(v\vert_{\partial\Omega}),v\vert_{\partial\Omega}>
\equiv \int_{\partial\Omega}(\Lambda_0-\Lambda_D)(v\vert_{\partial\Omega})(z)\,v(z)\,dS(z)
\tag {1.15}
$$
plays the central role, where $v$ is an arbitrary solution of the Helmholtz equation $\Delta v+k^2 v=0$ in $\Omega$.

The idea of the method of complementing function mentioned above suggests us the decomposition 
formulae for the form (1.15) stated below.

\proclaim{\noindent Theorem 1.2.} Let $v\in H^2(\Omega)$ satisfy
the Helmholtz equation $\Delta v+k^2v=0$ in $\Omega$.
We have
$$\begin{array}{ll}
\displaystyle
<(\Lambda_0-\Lambda_D)(v\vert_{\partial\Omega}),v\vert_{\partial\Omega}>
&
\displaystyle
=\Vert\nabla v\Vert_{L^2(D_n)}^2-k^2\Vert v\Vert_{L^2(D_n)}^2
\\
\\
\displaystyle
&
\displaystyle
\,\,\,
+\Vert \nabla(w+\epsilon_n)\Vert_{L^2(\Omega\setminus\overline{D})}^2-k^2\Vert (w+\epsilon_n)\Vert_{L^2(\Omega\setminus\overline{D})}^2
\\
\\
\displaystyle
&
\displaystyle
\,\,\,
-\Vert \nabla\epsilon_n\Vert_{L^2(\Omega\setminus\overline{D})}^2+k^2\Vert \epsilon_n\Vert_{L^2(\Omega\setminus\overline{D})}^2
\\
\\
\displaystyle
&
\displaystyle
\,\,\,
-\Vert\nabla v\Vert_{L^2(D_d)}^2+k^2\Vert v\Vert_{L^2(D_d)}^2
\end{array}
\tag {1.16}
$$
and
$$\begin{array}{ll}
\displaystyle
<(\Lambda_0-\Lambda_D)(v\vert_{\partial\Omega}),v\vert_{\partial\Omega}>
&
\displaystyle
=-\Vert\nabla v\Vert_{L^2(D_d)}^2
+k^2\Vert v\Vert_{L^2(D_d)}^2
\\
\\
\displaystyle
&
\displaystyle
\,\,\,
-\Vert \nabla(w+\epsilon_d)\Vert_{L^2(\Omega\setminus\overline{D})}^2
+k^2\Vert (w+\epsilon_d)\Vert_{L^2(\Omega\setminus\overline{D})}^2
\\
\\
\displaystyle
&
\displaystyle
\,\,\,
+\Vert \nabla\epsilon_d\Vert_{L^2(\Omega\setminus\overline{D})}^2
-k^2\Vert\epsilon_d\Vert_{L^2(\Omega\setminus\overline{D})}^2
\\
\\
\displaystyle
&
\displaystyle
\,\,\,
+\Vert\nabla v\Vert_{L^2(D_n)}^2-k^2\Vert v\Vert_{L^2(D_n)}^2,
\end{array}
\tag {1.17}
$$
where $w$, $\epsilon_n$ and $\epsilon_d$ are given by the solutions
$w_x$, $(\epsilon_x)_n$ and $(\epsilon_x)_d$ of  (1.6), (1.10) and (1.12) with $G(y,x)$ replaced by $v(y)$, respectively.

\endproclaim

\noindent
The formulae (1.16) and (1.17) are new and useful for establishing the probe method for Problem.
Those should be considered as the generalization of the well known decomposition formula in the case when $D_d=\emptyset$ or $D_n=\emptyset$, see \cite{Iwave} for the Helmholtz equation case. 
And also note that the expression of  the right-hand side on  (1.17) coincides with the one on (1.16) multiplied by $(-1)$ and replaced $(n,d)$ with $(d,n)$.

It should be emphasized that the proof of Theorem 1.2 given in Subsection 3.2 is independent of IPS.  Besides, the decomposition formulae (1.16) and (1.17) themselves would be valid also in the context of the weak solution.
However, formulae (1.16) and (1.17) should be considered as a byproduct of introducing the IPS function {\it at first}.
Without seeking the energy decomposition of IPS as done in \cite{IPS} one could never find the idea of the method of 
complementing functions to form (1.15).

This paper is organized as follows.  In Section 2 the proof of Theorem 1.1 is given.
The proof is based on integration by parts and clarifies the meaning of introducing the complementary
functions $(\epsilon_x)_n$ and $(\epsilon_x)_d$.
Section 3 is devoted to the integrated theory of the probe and singular sources methods.
In Subsection 3.1 a representation formula (3.2) of the indicator function for the probe method
as a limit of the indicator sequence (see Definition 3.1) in terms of the IPS function is established.   It is Theorem 3.1.  This together with Theorem 1.1
yields the Side A of the probe method which is concerned with blowing up of the indicator function 
on the surface of obstacles. 
In Subsection 3.2 first the proof of Theorem 1.2 together with its corollary is given.
Besides, In Subsection 3.3 it is shown that Theorem 1.2 yields the Side B
of the probe method which is concerned with blowing up of indicator sequence of the probe method
and stated as Theorem 3.2.
In Subsection 3.4 we will see that the singular sources method is included
in the IPS theory and in Subsection 3.4  it is shown that the singular sources method has the same side as the Side B of the probe method.  
Subsection 3.6 is devoted to a set of additional remarks related to the natural decomposition (1.5).
In the last section the conclusion and
some possible applications are briefly mentioned.

In Appendix we describe two lemmas which yield the upper estimate of the $L^2$-norm of the reflected solution $w_x$ and are used
in the proof of Corollary 1.1.

\section{Proof of Theorem 1.1}

First we show that the $w_x(y)$ which is the solution of  (1.6) has two expressions.
In what follows we always assume that $(x,y)\in\,(\Omega\setminus\overline{D})^2$.

\proclaim{\noindent Lemma 2.1.}  
It holds that 
$$\begin{array}{ll}
\displaystyle
w_x(y)
&
\displaystyle
=\int_{\partial\Omega}\frac{\partial}{\partial\nu}w_x(z)G(z,y)\,dS(z)\\
\\
\displaystyle
&
\displaystyle
\,\,\,
+\int_{D_n}\nabla G(z,x)\cdot\nabla G(z,y)\,dz-\int_{D_n}k^2 G(z,x)G(z,y)\,dz
\\
\\
\displaystyle
&
\displaystyle
\,\,\,
+\int_{\Omega\setminus\overline{D}}\nabla w_x(z)\cdot\nabla w_y(z)\,dz-\int_{\Omega\setminus\overline{D}}k^2w_x(z)w_y(z)\,dz
\\
\\
\displaystyle
&
\displaystyle
\,\,\,
-\int_{D_d}\nabla G(z,x)\cdot\nabla G(z,y)\,dz+\int_{D_d}k^2 G(z,x)G(z,y)\,dz
\\
\\
\displaystyle
&
\displaystyle
\,\,\,
+\int_{\partial D_d}\left(w_y(z)\frac{\partial}{\partial\nu}w_x(z)+w_x(z)\frac{\partial}{\partial\nu}w_y(z)\right)\,dS(z)
\end{array}
\tag {2.1}
$$
and
$$\begin{array}{ll}
\displaystyle
w_x(y)
&
\displaystyle
=\int_{\partial\Omega}\frac{\partial}{\partial\nu}w_x(z)G(z,y)\,dS(z)\\
\\
\displaystyle
&
\displaystyle
\,\,\,
-\int_{D_d}\nabla G(z,x)\cdot\nabla G(z,y)\,dz+\int_{D_d}k^2 G(z,x)G(z,y)\,dz
\\
\\
\displaystyle
&
\displaystyle
\,\,\,
-\int_{\Omega\setminus\overline{D}}\nabla w_x(z)\cdot\nabla w_y(z)\,dz
+\int_{\Omega\setminus\overline{D}}k^2w_x(z)w_y(z)\,dz
\\
\\
\displaystyle
&
\displaystyle
\,\,\,
+\int_{D_n}\nabla G(z,x)\cdot\nabla G(z,y)\,dz-\int_{D_n}k^2 G(z,x)G(z,y)\,dz
\\
\\
\displaystyle
&
\displaystyle
\,\,\,
-\int_{\partial D_n}\left(w_y(z)\frac{\partial}{\partial\nu}w_x(z)+w_x(z)\frac{\partial}{\partial\nu}w_y(z)\right)\,dS(z).
\end{array}
\tag {2.2}
$$
\endproclaim

{\it\noindent Proof.}
We start with the standard expression
$$\begin{array}{ll}
\displaystyle
w_x(y)
&
\displaystyle
=\int_{\partial\Omega}
\left(\frac{\partial}{\partial\nu}w_x(z)G(z,y)-
w_x(z)\frac{\partial}{\partial\nu}G(z,y)\right)\,dS(z)\\
\\
\displaystyle
&
\displaystyle
\,\,\,
-\int_{\partial D_n}
\left(\frac{\partial}{\partial\nu}w_x(z)G(z,y)-
w_x(z)\frac{\partial}{\partial\nu}G(z,y)\right)\,dS(z)
\\
\\
\displaystyle
&
\displaystyle
\,\,\,
-\int_{\partial D_d}
\left(\frac{\partial}{\partial\nu}w_x(z)G(z,y)-
w_x(z)\frac{\partial}{\partial\nu}G(z,y)\right)\,dS(z).
\end{array}
\tag {2.3}
$$
Applying the boundary conditions on $\partial\Omega$, $\partial D_d$ and $\partial D_n$ to (2.3), we obtain
$$\begin{array}{ll}
\displaystyle
w_x(y)
&
\displaystyle
=\int_{\partial\Omega}\frac{\partial}{\partial\nu}w_x(z)G(z,y)\,dS(z)\\
\\
\displaystyle
&
\displaystyle
\,\,\,
+\int_{\partial D_n}
\left(\frac{\partial}{\partial\nu}G(z,x)G(z,y)-
w_x(z)\frac{\partial}{\partial\nu}w_y(z)\right)\,dS(z)
\\
\\
\displaystyle
&
\displaystyle
\,\,\,
+\int_{\partial D_d}
\left(\frac{\partial}{\partial\nu}w_x(z)w_y(z)-G(z,x)\frac{\partial}{\partial\nu}G(z,y)\right)\,dS(z).
\end{array}
\tag {2.4}
$$
Since $x$ and $y$ outside both $D_d$ and $D_n$, we have, for $*=d,n$
$$\begin{array}{ll}
\displaystyle
\int_{\partial D_*}
\frac{\partial}{\partial\nu}G(z,x)G(z,y)\,dS(z)
&
\displaystyle
=\int_{D_*}\Delta G(z,x)G(z,y)\,dz+\int_{D_*}\nabla G(z,x)\cdot\nabla G(z,y)\,dz\\
\\
\displaystyle
&
\displaystyle
=\int_{D_*}\nabla G(z,x)\cdot\nabla G(z,y)\,dz-\int_{D_*}k^2 G(z,x)G(z,y)\,dz.
\end{array}
$$
Thus (2.4) becomes
$$\begin{array}{ll}
\displaystyle
w_x(y)
&
\displaystyle
=\int_{\partial\Omega}\frac{\partial}{\partial\nu}w_x(z)G(z,y)\,dS(z)\\
\\
\displaystyle
&
\displaystyle
\,\,\,
+\int_{D_n}\nabla G(z,x)\cdot\nabla G(z,y)\,dz-\int_{D_d}\nabla G(z,x)\cdot\nabla G(z,y)\,dz
\\
\\
\displaystyle
&
\displaystyle
\,\,\,
-\int_{D_n}k^2 G(z,x)G(z,y)\,dz+\int_{D_d}k^2 G(z,x)G(z,y)\,dz
\\
\\
\displaystyle
&
\displaystyle
\,\,\,
+\int_{\partial D_d}w_y(z)\frac{\partial}{\partial\nu}w_x(z)\,dS(z)
-\int_{\partial D_n}
w_x(z)\frac{\partial}{\partial\nu}w_y(z)\,dS(z).
\end{array}
\tag {2.5}
$$
Besides, we have
$$\begin{array}{l}
\displaystyle
\,\,\,\,\,\,
-\int_{\partial D_{n}}w_x(z)\frac{\partial}{\partial\nu}w_y(z)\,dS(z)-\int_{\partial D_{d}}w_x(z)\frac{\partial}{\partial\nu}w_y(z)\,dS(z)
\\
\\
\displaystyle
=\int_{\Omega\setminus\overline{D}}w_x(z)\Delta w_y(z)\,dz
+\int_{\Omega\setminus\overline{D}}\nabla w_x(z)\cdot\nabla w_y(z)\,dz
\\
\\
\displaystyle
=\int_{\Omega\setminus\overline{D}}\nabla w_x(z)\cdot\nabla w_y(z)\,dz-\int_{\Omega\setminus\overline{D}}k^2w_x(z)w_y(z)\,dz
\end{array}
\tag {2.6}
$$
and
$$\begin{array}{l}
\displaystyle
\,\,\,\,\,\,
-\int_{\partial D_{n}}w_y(z)\frac{\partial}{\partial\nu}w_x(z)\,dS(z)-\int_{\partial D_{d}}w_y(z)\frac{\partial}{\partial\nu}w_x(z)\,dS(z)
\\
\\
\displaystyle
=\int_{\Omega\setminus\overline{D}}\nabla w_x(z)\cdot\nabla w_y(z)\,dz-\int_{\Omega\setminus\overline{D}}k^2w_x(z)w_y(z)\,dz.
\end{array}
\tag {2.7}
$$
From (2.6) one has
$$\begin{array}{l}
\displaystyle
\,\,\,\,\,\,
-\int_{\partial D_{n}}w_x(z)\frac{\partial}{\partial\nu}w_y(z)\,dS(z)
\\
\\
\displaystyle
=\int_{\Omega\setminus\overline{D}}\nabla w_x(z)\cdot\nabla w_y(z)\,dz-\int_{\Omega\setminus\overline{D}}k^2w_x(z)w_y(z)\,dz
+\int_{\partial D_{d}}w_x(z)\frac{\partial}{\partial\nu}w_y(z)\,dS(z).
\end{array}
\tag {2.8}
$$
From (2.7) one has
$$\begin{array}{l}
\displaystyle
\,\,\,\,\,\,
\int_{\partial D_{d}}w_y(z)\frac{\partial}{\partial\nu}w_x(z)\,dS(z)
\\
\\
\displaystyle
=-\int_{\Omega\setminus\overline{D}}\nabla w_x(z)\cdot\nabla w_y(z)\,dz+\int_{\Omega\setminus\overline{D}}k^2w_x(z)w_y(z)\,dz
-\int_{\partial D_{n}}w_y(z)\frac{\partial}{\partial\nu}w_x(z)\,dS(z).
\end{array}
\tag {2.9}
$$
Thus one obtains the two representation of a single integral as follows.

From (2.8) we have
$$\begin{array}{l}
\displaystyle
\,\,\,\,\,\,
\int_{\partial D_d}w_y(z)\frac{\partial}{\partial\nu}w_x(z)\,dS(z)
-\int_{\partial D_n}
w_x(z)\frac{\partial}{\partial\nu}w_y(z)\,dS(z)
\\
\\
\displaystyle
=\int_{\partial D_d}\left(w_y(z)\frac{\partial}{\partial\nu}w_x(z)+w_x(z)\frac{\partial}{\partial\nu}w_y(z)\right)\,dS(z)
\\
\\
\displaystyle
\,\,\,
+
\int_{\Omega\setminus\overline{D}}\nabla w_x(z)\cdot\nabla w_y(z)\,dz-\int_{\Omega\setminus\overline{D}}k^2w_x(z)w_y(z)\,dz.
\end{array}
\tag {2.10}
$$
From (2.9) we have
$$\begin{array}{l}
\displaystyle
\,\,\,\,\,\,
\int_{\partial D_d}w_y(z)\frac{\partial}{\partial\nu}w_x(z)\,dS(z)
-\int_{\partial D_n}
w_x(z)\frac{\partial}{\partial\nu}w_y(z)\,dS(z)
\\
\\
\displaystyle
=-\int_{\partial D_n}\left(w_y(z)\frac{\partial}{\partial\nu}w_x(z)+w_x(z)\frac{\partial}{\partial\nu}w_y(z)\right)\,dS(z)
\\
\\
\displaystyle
\,\,\,
-
\int_{\Omega\setminus\overline{D}}\nabla w_x(z)\cdot\nabla w_y(z)\,dz+\int_{\Omega\setminus\overline{D}}k^2w_x(z)w_y(z)\,dz.
\end{array}
\tag {2.11}
$$
Subsituting  (2.10) and (2.11) into (2.5), we obtain (2.1) and (2.2).

\noindent
$\Box$

Next we show that the $w_x^1$ has the following expression.
\proclaim{\noindent Lemma 2.2.}  We have
$$\begin{array}{ll}
\displaystyle
w_x^1(y)
&
\displaystyle
=\int_{\Omega\setminus\overline{D}}\nabla w_x^1(z)\cdot\nabla w_y^1(z)\,dz-\int_{\Omega\setminus\overline{D}}k^2w_x^1(z)w_y^1(z)\,dz
\\
\\
\displaystyle
&
\displaystyle
\,\,\,
-\int_{\partial\Omega}
G(z,x)\frac{\partial}{\partial\nu}G(z,y)\,dS(z)
-\int_{\partial\Omega}G(z,x)\frac{\partial}{\partial\nu}w_y(z)\,dS(z).
\end{array}
\tag {2.12}
$$

\endproclaim

{\it\noindent Proof.}  To explain the reason for the introduction of the function $w_x^1$ step by step,
let us forget the set of  boundary conditions on  (1.7).

First same as (2.3) we start with the standard expression
$$\begin{array}{ll}
\displaystyle
w_x^1(y)
&
\displaystyle
=\int_{\partial\Omega}
\left(\frac{\partial}{\partial\nu}w_x^1(z)G(z,y)-
w_x^1(z)\frac{\partial}{\partial\nu}G(z,y)\right)\,dS(z)\\
\\
\displaystyle
&
\displaystyle
\,\,\,
-\int_{\partial D_n}
\left(\frac{\partial}{\partial\nu}w_x^1(z)G(z,y)-
w_x^1(z)\frac{\partial}{\partial\nu}G(z,y)\right)\,dS(z)
\\
\\
\displaystyle
&
\displaystyle
\,\,\,
-\int_{\partial D_d}
\left(\frac{\partial}{\partial\nu}w_x^1(z)G(z,y)-
w_x^1(z)\frac{\partial}{\partial\nu}G(z,y)\right)\,dS(z).
\end{array}
\tag {2.13}
$$
Here we impose the boundary conditions
$$
\left\{\begin{array}{ll}
\displaystyle
\frac{\partial}{\partial\nu}w_x^1(z)=0, & z\in\partial D_n,
\\
\\
\displaystyle
w_x^1(z)=G(z,x), & z\in\partial\Omega.
\end{array}
\right.
\tag {2.14}
$$
Then (2.13) becomes
$$\begin{array}{ll}
\displaystyle
w_x^1(y)
&
\displaystyle
=\int_{\partial\Omega}
\left(\frac{\partial}{\partial\nu}w_x^1(z)\,w_y^1(z)-
G(z,x)\frac{\partial}{\partial\nu}G(z,y)\right)\,dS(z)\\
\\
\displaystyle
&
\displaystyle
\,\,\,
+\int_{\partial D_n}
w_x^1(z)\frac{\partial}{\partial\nu}G(z,y)\,dS(z)
\\
\\
\displaystyle
&
\displaystyle
\,\,\,
-\int_{\partial D_d}
\left(\frac{\partial}{\partial\nu}w_x^1(z)\,G(z,y)-
w_x^1(z)\frac{\partial}{\partial\nu}G(z,y)\right)\,dS(z).
\end{array}
\tag {2.15}
$$
Here we have
$$\begin{array}{l}
\displaystyle
\,\,\,\,\,\,
\int_{\partial\Omega}
\frac{\partial}{\partial\nu}w_x^1(z)\,w_y^1(z)\,dS(z)
\\
\\
\displaystyle
=
\int_{\partial D_n}
\frac{\partial}{\partial\nu}w_x^1(z)\,w_y^1(z)\,dS(z)
+\int_{\partial D_d}
\frac{\partial}{\partial\nu}w_x^1(z)\,w_y^1(z)\,dS(z)
\\
\\
\displaystyle
\,\,\,
+\int_{\Omega\setminus\overline{D}}\Delta w_x^1(z)\,w_y^1(z)\,dz
\\
\\
\displaystyle
\,\,\,
+\int_{\Omega\setminus\overline{D}}\nabla w_x^1(z)\cdot\nabla w_y^1(z)\,dz
\\
\\
\displaystyle
=\int_{\Omega\setminus\overline{D}}\nabla w_x^1(z)\cdot\nabla w_y^1(z)\,dz
-\int_{\Omega\setminus\overline{D}}k^2w_x^1(z)w_y^1(z)\,dz
\\
\\
\displaystyle
\,\,\,
+\int_{\partial D_d}
\frac{\partial}{\partial\nu}w_x^1(z)\,w_y^1(z)\,dS(z).
\end{array}
$$
Thus (2.15) becomes
$$\begin{array}{ll}
\displaystyle
w_x^1(y)
&
\displaystyle
=\int_{\Omega\setminus\overline{D}}\nabla w_x^1(z)\cdot\nabla w_y^1(z)\,dz
-\int_{\Omega\setminus\overline{D}}k^2w_x^1(z)w_y^1(z)\,dz\\
\\
\displaystyle
&
\displaystyle
\,\,\,
-\int_{\partial\Omega}
G(z,x)\frac{\partial}{\partial\nu}G(z,y)\,dS(z)
+\int_{\partial D_n}
w_x^1(z)\frac{\partial}{\partial\nu}G(z,y)\,dS(z)
\\
\\
\displaystyle
&
\displaystyle
\,\,\,
+\int_{\partial D_d}
\frac{\partial}{\partial\nu}w_x^1(z)\,w_y^1(z)\,dS(z)
\\
\\
\displaystyle
&
\displaystyle
\,\,\,
-\int_{\partial D_d}
\left(\frac{\partial}{\partial\nu}w_x^1(z)\,G(z,y)-
w_x^1(z)\frac{\partial}{\partial\nu}G(z,y)\right)\,dS(z).
\end{array}
\tag {2.16}
$$
Here using the boundary condition of $w_y$ on $\partial D_n$, one has
$$\begin{array}{l}
\displaystyle
\,\,\,\,\,\,
\int_{\partial D_n}
w_x^1(z)\frac{\partial}{\partial\nu}G(z,y)\,dS(z)
\\
\\
\displaystyle
\displaystyle
=-\int_{\partial D_n}
w_x^1(z)\frac{\partial}{\partial\nu}w_y(z)\,dS(z)
\\
\\
\displaystyle
=\int_{\Omega\setminus\overline{D}} w_x^1(z)\Delta w_y(z)\,dz
-\int_{\partial\Omega}w_x^1(z)\frac{\partial}{\partial\nu}w_y(z)\,dS(z)
\\
\\
\displaystyle
\,\,\,
+\int_{\partial D_d}
w_x^1(z)\frac{\partial}{\partial\nu}w_y(z)\,dS(z)
\\
\\
\displaystyle
\,\,\,
+\int_{\Omega\setminus\overline{D}}\nabla w_x^1(z)\cdot\nabla w_y(z)\,dz\\
\\
\displaystyle
=-\int_{\partial\Omega}G(z,x)\frac{\partial}{\partial\nu}w_y(z)\,dS(z)
\\
\\
\displaystyle
\,\,\,
+\int_{\partial D_d}
w_x^1(z)\frac{\partial}{\partial\nu}w_y(z)\,dS(z)
+\int_{\Omega\setminus\overline{D}}\nabla w_x^1(z)\cdot\nabla w_y(z)\,dz
-\int_{\Omega\setminus\overline{D}}k^2w_x^1(z)w_y(z)\,dz.
\end{array}
$$
Besides we have
$$\begin{array}{ll}
\displaystyle
\int_{\Omega\setminus\overline{D}}\nabla w_x^1(z)\cdot\nabla w_y(z)\,dz
&
\displaystyle
=-\int_{\partial D_d}\frac{\partial}{\partial\nu}w_x^1(z)\,w_y(z)\,dS(z)-\int_{\Omega\setminus\overline{D}}\Delta w_x^1(z)w_y(z)\,dz
\\
\\
\displaystyle
&
\displaystyle
=\int_{\partial D_d}\frac{\partial}{\partial\nu}w_x^1(z)\,G(z,y)\,dS(z)
+\int_{\Omega\setminus\overline{D}}k^2w_x^1(z)w_y(z)\,dz.
\end{array}
$$
That is
$$
\displaystyle
\int_{\Omega\setminus\overline{D}}\nabla w_x^1(z)\cdot\nabla w_y(z)\,dz-\int_{\Omega\setminus\overline{D}}k^2w_x^1(z)w_y(z)\,dz
=\int_{\partial D_d}\frac{\partial}{\partial\nu}w_x^1(z)\,G(z,y)\,dS(z).
\tag {2.17}
$$
Note that we made use of the first boundary condition on $\partial D_n$ of (2.14) and 
the boundary condition for $w_y$ on $\partial D_d$ of (1.6).
Thus one gets
$$\begin{array}{l}
\displaystyle
\,\,\,\,\,\,
\int_{\partial D_n}
w_x^1(z)\frac{\partial}{\partial\nu}G(z,y)\,dS(z)
\\
\\
\displaystyle
=-\int_{\partial\Omega}G(z,x)\frac{\partial}{\partial\nu}w_y(z)\,dS(z)\\
\\
\displaystyle
\,\,\,
+\int_{\partial D_d}
w_x^1(z)\frac{\partial}{\partial\nu}w_y(z)\,dS(z)
+\int_{\partial D_d}\frac{\partial}{\partial\nu}w_x^1(z)\,G(z,y)\,dS(z).
\end{array}
$$
Therefore (2.16) becomes
$$\begin{array}{l}
\displaystyle
\,\,\,\,\,\,
w_x^1(y)
\\
\\
\displaystyle
=\int_{\Omega\setminus\overline{D}}\nabla w_x^1(z)\cdot\nabla w_y^1(z)\,dz-\int_{\Omega\setminus\overline{D}}k^2w_x^1(z)w_y^1(z)\,dz
\\
\\
\displaystyle
\,\,\,
-\int_{\partial\Omega}
G(z,x)\frac{\partial}{\partial\nu}G(z,y)\,dS(z)
-\int_{\partial\Omega}G(z,x)\frac{\partial}{\partial\nu}w_y(z)\,dS(z)\\
\\
\displaystyle
\,\,\,
+\int_{\partial D_d}
w_x^1(z)\frac{\partial}{\partial\nu}w_y(z)\,dS(z)
+\int_{\partial D_d}\frac{\partial}{\partial\nu}w_x^1(z)\,G(z,y)\,dS(z)
\\
\\
\displaystyle
\,\,\,
+\int_{\partial D_d}
\frac{\partial}{\partial\nu}w_x^1(z)\,w_y^1(z)\,dS(z)
\\
\\
\displaystyle
\,\,\,
-\int_{\partial D_d}
\left(\frac{\partial}{\partial\nu}w_x^1(z)\,G(z,y)-
w_x^1(z)\frac{\partial}{\partial\nu}G(z,y)\right)\,dS(z)
\\
\\
\displaystyle
=\int_{\Omega\setminus\overline{D}}\nabla w_x^1(z)\cdot\nabla w_y^1(z)\,dz-\int_{\Omega\setminus\overline{D}}k^2w_x^1(z)w_y^1(z)\,dz
\\
\\
\displaystyle
\,\,\,
-\int_{\partial\Omega}
G(z,x)\frac{\partial}{\partial\nu}G(z,y)\,dS(z)
-\int_{\partial\Omega}G(z,x)\frac{\partial}{\partial\nu}w_y(z)\,dS(z)\\
\\
\displaystyle
\,\,\,
+\int_{\partial D_d}
w_x^1(z)\frac{\partial}{\partial\nu}w_y(z)\,dS(z)
\\
\\
\displaystyle
\,\,\,
+\int_{\partial D_d}
\frac{\partial}{\partial\nu}w_x^1(z)\,w_y^1(z)\,dS(z)
+\int_{\partial D_d}
w_x^1(z)\frac{\partial}{\partial\nu}G(z,y)\,dS(z).
\end{array}
\tag {2.18}
$$
Here we impose the boundary condition of  $w_x^1$ and $w_y^1$ on $\partial D_d$:
$$
\begin{array}{ll}
\displaystyle
w_x^1(z)=w_y^1(z)=0, & z\in\partial D_d.
\end{array}
\tag {2.19}
$$
Thus (2.18) yields (2.12).

\noindent
$\Box$

\noindent
Note that the set of boundary conditions (2.14) and (2.19) coincides with that of (1.7).

From (2.1), (2.2) and (2.12) we immediately obtain the following two expressions for $W_x(y)$.

\proclaim{\noindent Proposition  2.1.}
It holds that
$$\begin{array}{l}
\displaystyle
\,\,\,\,\,\,
W_x(y)
\\
\\
\displaystyle
=\int_{D_n}\nabla G(z,x)\cdot\nabla G(z,y)\,dz-\int_{D_n}k^2 G(z,x)G(z,y)\,dz
\\
\\
\displaystyle
\,\,\,
+\int_{\Omega\setminus\overline{D}}\,\nabla w_x(z)\cdot\nabla w_y(z)\,dz-\int_{\Omega\setminus\overline{D}}k^2w_x(z)w_y(z)\,dz
\\
\\
\displaystyle
\,\,\,
-\int_{\partial\Omega}G(z,x)\frac{\partial}{\partial\nu}G(z,y)dS(z)
\\
\\
\displaystyle
\,\,\,
+\int_{\Omega\setminus\overline{D}}\nabla w_x^1(z)\cdot\nabla w_y^1(z)\,dz-\int_{\Omega\setminus\overline{D}}k^2w_x^1(z)w_y^1(z)\,dz
\\
\\
\displaystyle
\,\,\,
+\int_{\partial\Omega}G(z,y)\frac{\partial}{\partial\nu}w_x(z)\,dS(z)
-\int_{\partial\Omega}G(z,x)\frac{\partial}{\partial\nu}w_y(z)\,dS(z)
\\
\\
\displaystyle
\,\,\,
-\int_{D_d}\nabla G(z,x)\cdot\nabla G(z,y)\,dz+\int_{D_d}k^2G(z,x)G(z,y)\,dz
\\
\\
\displaystyle
\,\,\,
-\int_{\partial D_d}\left(G(z,y)\frac{\partial}{\partial\nu}w_x(z)+G(z,x)\frac{\partial}{\partial\nu}w_y(z)\right)\,dS(z)
\end{array}
\tag {2.20}
$$
and
$$\begin{array}{l}
\displaystyle
\,\,\,\,\,\,
W_x(y)
\\
\\
\displaystyle
=-\int_{D_d}\nabla G(z,x)\cdot\nabla G(z,y)\,dz+\int_{D_d}k^2G(z,x)G(z,y)\,dz
\\
\\
\displaystyle
\,\,\,-\int_{\Omega\setminus\overline{D}}\nabla w_x(z)\cdot\nabla w_y(z)\,dz
+\int_{\Omega\setminus\overline{D}}k^2w_x(z)w_y(z)\,dz
\\
\\
\displaystyle
\,\,\,
-\int_{\partial\Omega}G(z,x)\frac{\partial}{\partial\nu}G(z,y)dS(z)
\\
\\
\displaystyle
\,\,\,
+\int_{\Omega\setminus\overline{D}}\nabla w_x^1(z)\cdot\nabla w_y^1(z)\,dz-\int_{\Omega\setminus\overline{D}}k^2w_x^1(z)w_y^1(z)\,dz
\\
\\
\displaystyle
\,\,\,
+
\int_{\partial\Omega}
G(z,y)\frac{\partial}{\partial\nu}w_x(z)\,dS(z)
-\int_{\partial\Omega}G(z,x)\frac{\partial}{\partial\nu}w_y(z)\,dS(z)
\\
\\
\displaystyle
\,\,\,
+\int_{D_n}\nabla G(z,x)\cdot\nabla G(z,y)\,dz-\int_{D_n}k^2G(z,x)G(z,y)\,dz
\\
\\
\displaystyle
\,\,\,
+\int_{\partial D_n}\left(w_y(z)\frac{\partial}{\partial\nu}G(z,x)+w_x(z)\frac{\partial}{\partial\nu}G(z,y)\right)\,dS(z).
\end{array}
\tag {2.21}
$$
\endproclaim

Now let us explain the role of introducing the solutions of (1.10) and (1.12).
It is concerned with the last terms on  (2.20) and (2.21).
We call this technique the {\it method of complementing function}.

\proclaim{\noindent Lemma 2.3.}  We have
$$\begin{array}{l}
\displaystyle
\,\,\,\,\,\,
\int_{\partial D_d}G(z,y)\frac{\partial}{\partial\nu}w_x(z)\,dS(z)
+\int_{\partial D_d}G(z,x)\frac{\partial}{\partial\nu}w_y(z)\,dS(z)
\\
\\
\displaystyle
=-\int_{\Omega\setminus\overline{D}}
\left(\nabla(\epsilon_y)_n(z)\cdot\nabla w_x(z)+
\nabla(\epsilon_x)_n(z)\cdot\nabla w_y(z)\right)\,dz
\\
\\
\displaystyle
\,\,\,
+\int_{\Omega\setminus\overline{D}}
k^2\left((\epsilon_y)_n(z) w_x(z)+
(\epsilon_x)_n(z)w_y(z)\right)\,dz
\end{array}
\tag {2.22}
$$
and
$$\begin{array}{l}
\displaystyle
\,\,\,\,\,\,
\int_{\partial D_n}w_y(z)\frac{\partial}{\partial\nu}G(z,x)\,dS(z)
+\int_{\partial D_n}w_x(z)\frac{\partial}{\partial\nu}G(z,y)\,dS(z)
\\
\\
\displaystyle
=-\int_{\Omega\setminus\overline{D}}
\left(\nabla w_y(z)\cdot\nabla(\epsilon_x)_d(z)+
\nabla w_x(z)\cdot\nabla(\epsilon_y)_d(z)\right)\,dz\\
\\
\displaystyle
\,\,\,
+\int_{\Omega\setminus\overline{D}}
k^2\left(w_y(z) (\epsilon_x)_d(z)+
w_x(z) (\epsilon_y)_d(z)\right)\,dz.
\end{array}
\tag {2.23}
$$

\endproclaim

{\it\noindent Proof.}
First we rewrite the integral
$$\displaystyle
\int_{\partial D_d}G(z,y)\frac{\partial}{\partial\nu}w_x(z)\,dS(z).
$$
Using the equation (1.10), we have
$$\begin{array}{l}
\displaystyle
\,\,\,\,\,\,
\int_{\partial D_d}G(z,y)\frac{\partial}{\partial\nu}w_x(z)\,dS(z)
\\
\\
\displaystyle
=\int_{\partial D_d}(\epsilon_y)_n(z)\frac{\partial}{\partial\nu}w_x(z)\,dS(z)
+\int_{\partial D_n}(\epsilon_y)_n(z)\frac{\partial}{\partial\nu}w_x(z)\,dS(z)
-\int_{\partial\Omega}(\epsilon_y)_n(z)\frac{\partial}{\partial\nu}w_x(z)\,dS(z)
\\
\\
\displaystyle
=-\int_{\Omega\setminus\overline{D}}
(\epsilon_y)_n\Delta w_x(z)\,dz
-\int_{\Omega\setminus\overline{D}}
\nabla(\epsilon_y)_n\cdot\nabla w_x(z)\,dz
\\
\\
\displaystyle
=-\int_{\Omega\setminus\overline{D}}
\nabla(\epsilon_y)_n\cdot\nabla w_x(z)\,dz
+\int_{\Omega\setminus\overline{D}}k^2(\epsilon_y)_n(z)w_x(z)\,dz.
\end{array}
$$
Interchanging $x$ and $y$, we obtain the expression (2.22).

Second we rewrite the integral
$$\displaystyle
\int_{\partial D_n}w_x(z)\frac{\partial}{\partial\nu}G(z,x)\,dS(z).
$$
Using the equation (1.12), we have
$$\begin{array}{l}
\displaystyle
\,\,\,\,\,\,
\int_{\partial D_n}w_x(z)\frac{\partial}{\partial\nu}G(z,y)\,dS(z)
\\
\\
\displaystyle
=\int_{\partial D_n}w_x(z)\frac{\partial}{\partial\nu}(\epsilon_y)_d(z)\,dS(z)
+\int_{\partial D_d}w_x(z)\frac{\partial}{\partial\nu}(\epsilon_y)_d(z)\,dS(z)
-\int_{\partial\Omega}w_x(z)\frac{\partial}{\partial\nu}(\epsilon_y)_d(z)\,dS(z)
\\
\\
\displaystyle
=-\int_{\Omega\setminus\overline{D}}
w_x\Delta(\epsilon_y)_d(z)\,dz
-\int_{\Omega\setminus\overline{D}}
\nabla w_x(z)\cdot\nabla(\epsilon_y)_d(z)\,dz
\\
\\
\displaystyle
=-\int_{\Omega\setminus\overline{D}}
\nabla w_x(z)\cdot\nabla(\epsilon_y)_d(z)\,dz+\int_{\Omega\setminus\overline{D}}k^2 w_x(z)(\epsilon_y)_d(z)\,dz.
\end{array}
$$
Interchanging $x$ and $y$, we obtain the expression (2.23).

\noindent
$\Box$

Thus (2.20) together with (2.22) yields
$$\begin{array}{l}
\displaystyle
\,\,\,\,\,\,
W_x(y)
\\
\\
\displaystyle
=\int_{D_n}\nabla G(z,x)\cdot\nabla G(z,y)\,dz-\int_{D_n}k^2G(z,x)G(z,y)\,dz\\
\\
\displaystyle
\,\,\,
+\int_{\Omega\setminus\overline{D}}\,\nabla w_x(z)\cdot\nabla w_y(z)\,dz-\int_{\Omega\setminus\overline{D}}k^2 w_x(z)w_y(z)\,dz
\\
\\
\displaystyle
\,\,\,
-\int_{\partial\Omega}\,G(z,x)\frac{\partial}{\partial\nu}G(z,y)\,dS(z)\\
\\
\displaystyle
\,\,\,
+\int_{\Omega\setminus\overline{D}}\nabla w_x^1(z)\cdot\nabla w_y^1(z)\,dz
-\int_{\Omega\setminus\overline{D}}k^2 w_x^1(z)w_y^1(z)\,dz
\\
\\
\displaystyle
\,\,\,
+\int_{\partial\Omega}G(z,y)\frac{\partial}{\partial\nu}w_x(z)\,dS(z)
-\int_{\partial\Omega}G(z,x)\frac{\partial}{\partial\nu}w_y(z)\,dS(z)
\\
\\
\displaystyle
\,\,\,
-\int_{D_d}\nabla G(z,x)\cdot\nabla G(z,y)\,dz
+\int_{D_d}k^2G(z,x)G(z,y)\,dz
\\
\\
\displaystyle
\,\,\,
+\int_{\Omega\setminus\overline{D}}
\left(\nabla(\epsilon_y)_n(z)\cdot\nabla w_x(z)+
\nabla(\epsilon_x)_n(z)\cdot\nabla w_y(z)\right)\,dz
\\
\\
\displaystyle
-\int_{\Omega\setminus\overline{D}}
k^2\left((\epsilon_y)_n(z) w_x(z)+
(\epsilon_x)_n(z) w_y(z)\right)\,dz.
\end{array}
\tag {2.24}
$$
And (2.21) together with (2.23) yields
$$\begin{array}{l}
\displaystyle
\,\,\,\,\,\,
W_x(y)
\\
\\
\displaystyle
=-\int_{D_d}\nabla G(z,x)\cdot\nabla G(z,y)\,dz
+\int_{D_d}k^2G(z,x)G(z,y)\,dz
\\
\\
\displaystyle
\,\,\,
-\int_{\Omega\setminus\overline{D}}\nabla w_x(z)\cdot\nabla w_y(z)\,dz
+\int_{\Omega\setminus\overline{D}}k^2w_x(z) w_y(z)\,dz
\\
\\
\displaystyle
\,\,\,
-\int_{\partial\Omega}\,G(z,x)\frac{\partial}{\partial\nu}G(z,y)\,dS(z)
\\
\\
\displaystyle
\,\,\,
+\int_{\Omega\setminus\overline{D}}\nabla w_x^1(z)\cdot\nabla w_y^1(z)\,dz
-\int_{\Omega\setminus\overline{D}}k^2w_x^1(z) w_y^1(z)\,dz
\\
\\
\displaystyle
\,\,\,
+
\int_{\partial\Omega}G(z,y)\frac{\partial}{\partial\nu}w_x(z)\,dS(z)
-\int_{\partial\Omega}G(z,x)\frac{\partial}{\partial\nu}w_y(z)\,dS(z)
\\
\\
\displaystyle
\,\,\,
+\int_{D_n}\nabla G(z,x)\cdot\nabla G(z,y)\,dz-\int_{D_n}k^2G(z,x)G(z,y)\,dz
\\
\\
\displaystyle
\,\,\,
-\int_{\Omega\setminus\overline{D}}
\left(\nabla w_y(z)\cdot\nabla(\epsilon_x)_d(z)+
\nabla w_x(z)\cdot\nabla(\epsilon_y)_d(z)\right)\,dz
\\
\\
\displaystyle
\,\,\,
+\int_{\Omega\setminus\overline{D}}k^2
\left(w_y(z) (\epsilon_x)_d(z)+
w_x(z) (\epsilon_y)_d(z)\right)\,dz.
\end{array}
\tag {2.25}
$$

Letting $x=y$ in (2.24) and (2.25), we obtain
$$\begin{array}{l}
\displaystyle
\,\,\,\,\,\,
W_x(x)
\\
\\
\displaystyle
=\Vert\nabla G(\,\cdot\,,x)\Vert_{L^2(D_n)}^2-k^2\Vert G(\,\cdot\,,x)\Vert_{L^2(D_n)}^2
+\Vert \nabla w_x\Vert_{L^2(\Omega\setminus\overline{D})}^2-k^2\Vert w_x\Vert_{L^2(\Omega\setminus\overline{D})}^2
\\
\\
\displaystyle
\,\,\,
-\int_{\partial\Omega}\,G(z,x)\frac{\partial}{\partial\nu}G(z,x)\,dS(z)
+
\Vert \nabla w_x^1\Vert_{L^2(\Omega\setminus\overline{D})}^2-k^2\Vert w_x^1\Vert_{L^2(\Omega\setminus\overline{D})}^2
\\
\\
\displaystyle
\,\,\,
-\Vert\nabla G(\,\cdot\,,x)\Vert_{L^2(D_d)}^2+k^2\Vert G(\,\cdot\,,x)\Vert_{L^2(D_d)}^2
\\
\\
\displaystyle
\,\,\,
+2\int_{\Omega\setminus\overline{D}}
\nabla(\epsilon_x)_n(z)\cdot\nabla w_x(z)\,dz
-2\int_{\Omega\setminus\overline{D}}
k^2w_x(z) (\epsilon_x)_n(z)\,dz
\end{array}
\tag {2.26}
$$
and
$$\begin{array}{l}
\displaystyle
\,\,\,\,\,\,
W_x(x)
\\
\\
\displaystyle
=-\Vert\nabla G(\,\cdot\,,x)\Vert_{L^2(D_d)}^2+k^2\Vert G(\,\cdot\,,x)\Vert_{L^2(D_d)}^2
-
\Vert \nabla w_x\Vert_{L^2(\Omega\setminus\overline{D})}^2
+k^2\Vert w_x\Vert_{L^2(\Omega\setminus\overline{D})}^2
\\
\\
\displaystyle
\,\,\,
-\int_{\partial\Omega}\,G(z,x)\frac{\partial}{\partial\nu}G(z,x)\,dS(z)
+
\Vert \nabla w_x^1\Vert_{L^2(\Omega\setminus\overline{D})}^2
-k^2\Vert \nabla w_x^1\Vert_{L^2(\Omega\setminus\overline{D})}^2
\\
\\
\displaystyle
\,\,\,
+\Vert\nabla G(\,\cdot\,,x)\Vert_{L^2(D_n)}^2
-k^2\Vert G(\,\cdot\,,x)\Vert_{L^2(D_n)}^2
\\
\\
\displaystyle
\,\,\,
-2\int_{\Omega\setminus\overline{D}}
\nabla w_x(z)\cdot\nabla(\epsilon_x)_d(z)\,dz
+2\int_{\Omega\setminus\overline{D}}
k^2w_x(z) (\epsilon_x)_d(z)\,dz.
\end{array}
\tag {2.27}
$$

Here rewrite
$$\begin{array}{l}
\displaystyle
\,\,\,\,\,\,
\Vert \nabla w_x\Vert_{L^2(\Omega\setminus\overline{D})}^2
+2\int_{\Omega\setminus\overline{D}}
\nabla(\epsilon_x)_n\cdot\nabla w_x(z)\,dz
\\
\\
\displaystyle
=\Vert \nabla(w_x+(\epsilon_x)_n)\Vert_{L^2(\Omega\setminus\overline{D})}^2
-\Vert \nabla(\epsilon_x)_n\Vert_{L^2(\Omega\setminus\overline{D})}^2
\end{array}
\tag {2.28}
$$
and
$$\begin{array}{l}
\displaystyle
\,\,\,\,\,\,
\Vert \nabla w_x\Vert_{L^2(\Omega\setminus\overline{D})}^2
+2\int_{\Omega\setminus\overline{D}}
\nabla w_x\cdot\nabla(\epsilon_x)_d(z)\,dz
\\
\\
\displaystyle
=\Vert \nabla(w_x+(\epsilon_x)_d)\Vert_{L^2(\Omega\setminus\overline{D})}^2
-\Vert \nabla(\epsilon_x)_d\Vert_{L^2(\Omega\setminus\overline{D})}^2.
\end{array}
\tag {2.29}
$$
Rewrite also as
$$\begin{array}{l}
\displaystyle
\,\,\,\,\,\,
\Vert w_x\Vert_{L^2(\Omega\setminus\overline{D})}^2
+2\int_{\Omega\setminus\overline{D}}
(\epsilon_x)_n w_x(z)\,dz
\\
\\
\displaystyle
=\Vert w_x+(\epsilon_x)_n\Vert_{L^2(\Omega\setminus\overline{D})}^2
-\Vert (\epsilon_x)_n\Vert_{L^2(\Omega\setminus\overline{D})}^2
\end{array}
\tag {2.30}
$$
and
$$\begin{array}{l}
\displaystyle
\,\,\,\,\,\,
\Vert w_x\Vert_{L^2(\Omega\setminus\overline{D})}^2
+2\int_{\Omega\setminus\overline{D}}
w_x (\epsilon_x)_d(z)\,dz
\\
\\
\displaystyle
=\Vert w_x+(\epsilon_x)_d\Vert_{L^2(\Omega\setminus\overline{D})}^2
-\Vert (\epsilon_x)_d\Vert_{L^2(\Omega\setminus\overline{D})}^2.
\end{array}
\tag {2.31}
$$
Then from  (2.26), (2.27), (2.28), (2.29), (2.30) and (2.31) we obtain (1.9) and (1.11) of Theorem 1.1.

$\quad$

{\bf\noindent Remark 2.1.}
Assume that $D_d=\emptyset$.  This the purely Neumann obstacle case.  Then the equation (2.17) becomes
$$\displaystyle
\int_{\Omega\setminus\overline{D}}\nabla w_x^1(z)\cdot\nabla w_y(z)\,dz-\int_{\Omega\setminus\overline{D}}k^2w_x^1(z)w_y(z)\,dz=0.
\tag {2.32}
$$
It is easy to see that equation (2.32) combined with (i) of Corollary 1.2 makes the representation of IPS function so simple:
$$\begin{array}{ll}
\displaystyle
W_x(x)
&
\displaystyle
=\Vert\nabla W_x\Vert_{L^2(\Omega\setminus\overline{D}_n)}^2-k^2\Vert W_x\Vert_{L^2(\Omega\setminus\overline{D_n})}^2
\\
\\
\displaystyle
&
\displaystyle
\,\,\,
+\Vert\nabla G(\,\cdot\,-x)\Vert_{L^2(D_n)}^2-k^2\Vert G(\,\cdot\,-x)\Vert_{L^2(D_n)}^2
-\int_{\partial\Omega}\,G(z-x)\frac{\partial}{\partial\nu}G(z-x)\,dS(z).
\end{array}
$$
This is an extension of  the expression of IPS function given in Remark 1.7 of \cite{IPS} to the case $k\not=0$.
It seems, in the case $D_d\not=\emptyset$ one cannot expect such a simple expression.

$\quad$

\section{Integrated theory}

\subsection{IPS to Side A of probe method}

In this subsection we derive the probe method via the integrated theory of the probe and singular sources methods.
We fix an arbitrary ${\cal G}=\{G(\,\cdot\,,x)\}_{x\in\Omega}$ given by (1.2) unless otherwise specified.

First we recall the notion of a needle.  Given $x\in\Omega$ let $N_x$ denote the set of all non-self intersecting piecewise linear curves $\sigma$ connecting a point on $\partial\Omega$
and $x$ such that other points on $\sigma$ are in $\Omega$.  We call each member in $N_x$ a needle with a tip at $x$.

$\quad$

{\bf\noindent Definition 3.1.}
Given $x\in\Omega$ and $\sigma\in N_x$
a sequence $\{v_n\}$ of $H^2(\Omega)$ functions is called a needle sequence for $(x,\sigma)$ {\it based on} ${\cal G}$
if each $v_n$ satisfies the Helmholtz equation $\Delta v+k^2 v=0$ in $\Omega$
and $\{v_n\}$ converges to $G(\,\cdot\,,x)$ in $H^2_{\text{loc}}(\Omega\setminus\sigma)$.
Then the sequence given by
$$\displaystyle
<(\Lambda_0-\Lambda_D)(v_n\vert_{\partial\Omega}),v_n\vert_{\partial\Omega}>
\equiv
\int_{\partial\Omega}(\Lambda_0-\Lambda_D)(v_n\vert_{\partial\Omega})(z)v_n(z)\,dS(z),
$$
is called the {\it indicator sequence} for the probe method.

$\quad$

\noindent
Note that, hereafter, unless otherwise specified we assume that Assumption 1 for the case $D=\emptyset$ is also satisfied.  
This ensures not only the well-definedness of $\Lambda_0$ but also the existence
of the needle sequence for an arbitrary needle \cite{Iwave}.

The Side A of the probe method starts with the convergence property
of the indicator sequence as described below.

\proclaim{\noindent Theorem 3.1.}  Let $x\in\Omega\setminus\overline{D}$ and
$\sigma\in N_x$.  Let $\{v_n\}$ be an arbitrary needle sequence for $(x,\sigma)$ based on ${\cal G}$.
If $\sigma\cap\overline{D}=\emptyset$, then we have
$$\displaystyle
\lim_{n\rightarrow\infty}<(\Lambda_0-\Lambda_D)(v_n\vert_{\partial\Omega}),v_n\vert_{\partial\Omega}>
=I(x),
\tag {3.1}
$$
where 
$$\displaystyle
I(x)=W_x(x)-<\Lambda_D(G(\,\cdot\,,x)\vert_{\partial\Omega}),G(\,\cdot\,,x)\vert_{\partial\Omega}>
+\int_{\partial\Omega}\,\frac{\partial}{\partial\nu}G(z,x)G(z,x)\,dS(z).
\tag {3.2}
$$

\endproclaim
{\it\noindent Proof.}
Fix $x\in\Omega\setminus\overline{D}$.
First we show that the limit of the left-hand side on (3.1) exists and its limit has the expression
$$
\displaystyle
\lim_{n\rightarrow\infty}<(\Lambda_0-\Lambda_D)(v_n\vert_{\partial\Omega},v_n\vert_{\partial\Omega}>
=w_x(x)-\int_{\partial\Omega}
\frac{\partial}{\partial\nu}w_x(z)G(z,x)\,dS(z).
\tag {3.3}
$$
Define
$$\begin{array}{ll}
\displaystyle
G_n(z,x)=G(z,x)-v_n(z),
&
\displaystyle
z\in\Omega.
\end{array}
\tag {3.4}
$$
The form of $G_n(\,\cdot\,,x)$ together with Green's theorem yields
an expression of $w_n=w(z)$ at $z=x$, which is the solution of 
$$
\left\{
\begin{array}{ll}
\displaystyle
\Delta w+k^2w=0, & z\in\Omega\setminus\overline{D},
\\
\\
\displaystyle
\frac{\partial w}{\partial\nu}=-\frac{\partial v_n}{\partial\nu}, & z\in\partial D_n,\\
\\
\displaystyle
w=-v_n, & z\in\partial D_d,
\\
\\
\displaystyle
w=0, & z\in\partial\Omega.
\end{array}
\right.
$$
That is
$$\begin{array}{ll}
\displaystyle
w_n(x)
&
\displaystyle
=\int_{\partial\Omega}
\frac{\partial}{\partial\nu}w_n(z)G_n(z,x)\,dS(z)\\
\\
\displaystyle
&
\displaystyle
-\int_{\partial D_n}
\left(\frac{\partial}{\partial\nu}w_n(z)G_n(z,x)-
w_n(z)\frac{\partial}{\partial\nu}G_n(z,x)\right)\,dS(z)
\\
\\
\displaystyle
&
\displaystyle
-\int_{\partial D_d}
\left(\frac{\partial}{\partial\nu}w_n(z)G_n(z,x)-
w_n(z)\frac{\partial}{\partial\nu}G_n(z,x)\right)\,dS(z).
\end{array}
\tag {3.5}
$$
By Definition 3.1 and the assumption $\sigma\cap\overline{D}=\emptyset$, we have, as $n\rightarrow\infty$
$G_n(\,\cdot\,,x)\rightarrow 0$ in $H^2(D)$
and the well-posedness, we have $w_n\rightarrow w_x$ in $H^2(\Omega\setminus\overline{D})$.
Therefore it follows from these and the Sobolev embedding, letting $n\rightarrow\infty$ of (3.5), we obtain
$$\begin{array}{ll}
\displaystyle
w_x(x)
&
\displaystyle
=\lim_{n\rightarrow\infty}\int_{\partial\Omega}
\frac{\partial}{\partial\nu}w_n(z)G_n(z,x)\,dS(z)
\\
\\
\displaystyle
&
\displaystyle
=\lim_{n\rightarrow\infty}
\left(\int_{\partial\Omega}\frac{\partial w_n}{\partial\nu\,}G(z,x)\,dS(z)
-
\int_{\partial\Omega}\frac{\partial w_n}{\partial\nu}\,v_n(z)\,dS(z)
\right).
\end{array}
\tag {3.6}
$$
Since $\frac{\partial w_n}{\partial\nu}\rightarrow\frac{\partial w_x}{\partial\nu}$ in $H^{\frac{1}{2}}(\partial\Omega)$,
the first term of the right-hand side on (3.6) is convergent and the limit is given by
$$\displaystyle
\int_{\partial\Omega}\frac{\partial w_x}{\partial\nu\,}G(z,x)\,dS(z).
$$
Therefore the second term of the right-hand side on (3.6) is also convergent and its limit
satisfies
$$\displaystyle
w_x(x)=\int_{\partial\Omega}\frac{\partial w_x}{\partial\nu\,}G(z,x)\,dS(z)
-\lim_{n\rightarrow\infty}\int_{\partial\Omega}\frac{\partial w_n}{\partial\nu}\,v_n(z)\,dS(z).
$$
Using this and the trivial expression
$$\begin{array}{ll}
\displaystyle
\frac{\partial}{\partial\nu}w_n(z)=-(\Lambda_0-\Lambda_D)(v_n\vert_{\partial\Omega})(z),
&
z\in\partial\Omega,
\end{array}
$$
we obtain
$$
\begin{array}{ll}
\displaystyle
w_x(x)
&
\displaystyle
=
\int_{\partial\Omega}\frac{\partial w_x}{\partial\nu}(z)\,
G(z,x)\,dS(z)
+\lim_{n\rightarrow\infty}<(\Lambda_0-\Lambda_D)(v_n\vert_{\partial\Omega}),v_n\vert_{\partial\Omega}>.
\end{array}
$$
This is nothing but (3.3).

Next we show that the right-hand side of the formula (3.3) coinsides with that of  formula (3.2).
Recalling  (2.12) of Lemma 2.2 with $x=y$, we have
$$\begin{array}{l}
\displaystyle
\,\,\,\,\,\,
-\int_{\partial\Omega}G(z,x)\frac{\partial}{\partial\nu}w_x(z)\,dS(z)\\
\\
\displaystyle
=w_x^1(x)-\Vert\nabla w_x^1\Vert_{L^2(\Omega\setminus\overline{D})}^2
+k^2\Vert w_x^1\Vert_{L^2(\Omega\setminus\overline{D})}^2
+\int_{\partial\Omega}
G(z,x)\frac{\partial}{\partial\nu}G(z,x)\,dS(z).
\end{array}
\tag {3.7}
$$
Besides from the equation (1.7) we have
$$\begin{array}{ll}
\displaystyle
<\Lambda_D(G(\,\cdot\,,x)\vert_{\partial\Omega},G(\,\cdot\,,x)\vert_{\partial\Omega}>
&
\displaystyle
=\int_{\partial\Omega}\frac{\partial}{\partial\nu}w_x^1(z)\,w_x^1(z)\,dS(z)
\\
\\
\displaystyle
&
\displaystyle
=\Vert\nabla w_x^1\Vert_{L^2(\Omega\setminus\overline{D})}^2
-k^2\Vert w_x^1\Vert_{L^2(\Omega\setminus\overline{D})}^2.
\end{array}
\tag {3.8}
$$
From (3.7) and (3.8) together with (1.8) we see that
$$\begin{array}{l}
\displaystyle
\,\,\,\,\,\,
\displaystyle
w_x(x)-\int_{\partial\Omega}\frac{\partial}{\partial\nu}w_x(z)G(z,x)\,dS(z)
\\
\\
\displaystyle
=w_x(x)+w_x^1(x)-\Vert\nabla w_x^1\Vert_{L^2(\Omega\setminus\overline{D})}^2
+k^2\Vert w_x^1\Vert_{L^2(\Omega\setminus\overline{D})}^2
+\int_{\partial\Omega}
G(z,x)\frac{\partial}{\partial\nu}G(z,x)\,dS(z)
\\
\\
\displaystyle
=W_x(x)-<\Lambda_D(G(\,\cdot\,,x)\vert_{\partial\Omega},G(\,\cdot\,,x)\vert_{\partial\Omega}>
+\int_{\partial\Omega}
G(z,x)\frac{\partial}{\partial\nu}G(z,x)\,dS(z).
\end{array}
\tag {3.9}
$$
A combination of  (3.3) and (3.9) yields the desired conclusion.

\noindent
$\Box$

$\quad$

{\bf\noindent Definition 3.2.}
The function $I(x)$ appeared as the limit (3.1) and expressed as (3.2)  is called the {\it indicator function for the probe method}
based on ${\cal G}$.

$\quad$

\noindent
The formula (3.1) should be understood as a computation formula of the indicator function by using $\Lambda_0-\Lambda_D$.
Besides, this shows that IPS function $W_x(x)$ can be calculated from $\Lambda_D$
(and $\Lambda_0$ which can be calculated in advance) acting on the needle sequences from the surface
$\partial\Omega$ to inside.

The second and third terms of the right-hand side on (3.2)
are bounded when $x$ is away from $\partial\Omega$, by virtue of (1.2) and (1.3).
Thus it follows from Corollary 1.1 and (3.2) that

$\quad$

\noindent
(i)  $\lim_{x\rightarrow a\in\partial D_n}I(x)=\infty$.

\noindent
(ii)  $\lim_{x\rightarrow b\in\partial D_d}I(x)=-\infty$.

$\quad$

Then it follows from Theorem 1.1, (3.2) and (3.8) that the $I(x)$ has two expressions.

$$\begin{array}{ll}
\displaystyle
I(x)
&
\displaystyle
=\Vert\nabla G(\,\cdot\,,x)\Vert_{L^2(D_n)}^2-k^2\Vert G(\,\cdot\,,x)\Vert_{L^2(D_n)}^2
\\
\\
\displaystyle
&
\displaystyle
\,\,\,
+\Vert \nabla(w_x+(\epsilon_x)_n)\Vert_{L^2(\Omega\setminus\overline{D})}^2
-k^2\Vert w_x+(\epsilon_x)_n\Vert_{L^2(\Omega\setminus\overline{D})}^2
\\
\\
\displaystyle
&
\displaystyle
\,\,\,
-\Vert \nabla(\epsilon_x)_n\Vert_{L^2(\Omega\setminus\overline{D})}^2
+k^2\Vert (\epsilon_x)_n\Vert_{L^2(\Omega\setminus\overline{D})}^2
\\
\\
\displaystyle
&
\displaystyle
\,\,\,
-\Vert\nabla G(\,\cdot\,,x)\Vert_{L^2(D_d)}^2+k^2\Vert G(\,\cdot\,,x)\Vert_{L^2(D_d)}^2,
\end{array}
\tag {3.10}
$$

$$\begin{array}{ll}
\displaystyle
I(x)
&
\displaystyle
=-\Vert\nabla G(\,\cdot\,,x)\Vert_{L^2(D_d)}^2+k^2\Vert G(\,\cdot\,,x)\Vert_{L^2(D_d)}^2
\\
\\
\displaystyle
&
\displaystyle
\,\,\,
-\Vert \nabla(w_x+(\epsilon_x)_d)\Vert_{L^2(\Omega\setminus\overline{D})}^2+
k^2\Vert w_x+(\epsilon_x)_d\Vert_{L^2(\Omega\setminus\overline{D})}^2
\\
\\
\displaystyle
&
\displaystyle
\,\,\,
+\Vert \nabla(\epsilon_x)_d\Vert_{L^2(\Omega\setminus\overline{D})}^2
-k^2\Vert (\epsilon_x)_d\Vert_{L^2(\Omega\setminus\overline{D})}^2
\\
\\
\displaystyle
&
\displaystyle
\,\,\,
+\Vert\nabla G(\,\cdot\,,x)\Vert_{L^2(D_n)}^2
-k^2\Vert G(\,\cdot\,,x)\Vert_{L^2(D_n)}^2
.\end{array}
\tag {3.11}
$$
Then we can easily check that, for each $\epsilon_i>0$, $i=1,2$
$$\displaystyle
\sup_{x\in\Omega\setminus\overline{D}, \text{dist}\,(x,\partial D)>\epsilon_1,
\,\text{dist}\,(x,\partial\Omega)>\epsilon_2}\,\vert I(x)\vert <\infty.
\tag {3.12}
$$

So the convergence of the indicator sequence (3.1),  blowing up property of the $I(x)$  mentioned  (i) and (ii) above
and (3.12) establish the Side A of the probe method for the mixed obstacle case.  

The point that should be emphasized
is:  from IPS we obtained that $I(x)$ as the limit of the indicator sequence takes the expressions 
(3.10) and (3.11).  Is sholud be also pointed out that,  the expression (3.11) coincides with $(-1)$-times the expression
(3.10) replaced with $(n,d)$ with $(d,n)$.

$\quad$

{\bf\noindent Remark 3.1.}
If ${\cal G}={\cal G}^0$, from the well-posedness of boundary value problems (1.6), (1.10) and (1.12)
one can relax  (3.12) as: for each $\epsilon>0$
$$\displaystyle
\sup_{x\in\Omega\setminus\overline{D}, \text{dist}\,(x,\partial D)>\epsilon}\,\vert I(x)\vert <\infty.
$$

$\quad$

{\bf\noindent Definition 3.3.}
For an arbitrary point $x\in\Omega\setminus\overline{D}$ define
$$
\displaystyle
I^1(x)=<\Lambda_D(G(\,\cdot\,,x)\vert_{\partial\Omega}),G(\,\cdot\,,x)\vert_{\partial\Omega}>-\int_{\partial\Omega}\,\frac{\partial}{\partial\nu}G(z,x)G(z,x)\,dS(z).
\tag {3.13}
$$

$\quad$

\noindent
In principle, it is possible to calculate $I^1(x)$ in advance from given $\Lambda_D$ without probing.
Besides, if ${\cal G}={\cal G}^0$, it follows from (3.8) we have the energy
integral expression of  (3.13):
$$
\displaystyle
I^1(x)=\Vert\nabla w_x^1\Vert_{L^2(\Omega\setminus\overline{D})}^2-k^2\Vert w_x^1\Vert_{L^2(\Omega\setminus\overline{D})}^2-\int_{\partial\Omega}\,\frac{\partial}{\partial\nu}G(z,x)G(z,x)\,dS(z).
$$
In particular, if $k=0$, then by (1.14) this becomes
$$\begin{array}{ll}
\displaystyle
I^1(x)=\Vert\nabla w_x^1\Vert_{L^2(\Omega\setminus\overline{D})}^2
+\Vert\nabla G(\,\cdot\,-x)\Vert_{L^2(\Bbb R^3\setminus\overline{\Omega})}^2, & x\in\Omega\setminus\overline{D}.
\end{array}
$$

Finally by (3.2) and (3.13) we have the {\it inner decomposition} of IPS function:
$$\begin{array}{ll}
\displaystyle
W_x(x)=I(x)+I^1(x), & x\in\Omega\setminus\overline{D}.
\end{array}
\tag {3.14}
$$
The equations (1.8) and (3.14) give us two ways of decomposition of IPS function.

$\quad$

{\bf\noindent Remark 3.2.}
The two types of the decompositions (3.10) and (3.11) suggest the replacement:
$$\displaystyle
G(\,\cdot\,,x)\rightarrow v,
$$
where $v$ is an arbitrary solution of the Helmholtz equation in $\Omega$.
Theorem 1.2 can be considered as an example of the validity of this replacement.

\subsection{Proof of Theorem 1.2 and a corollary}

{\it\noindent Proof of Theorem 1.2.}
First we prove the validity of  (1.16).
We have
$$\begin{array}{l}
\displaystyle
\,\,\,\,\,\,
\Vert\nabla(w+\epsilon_n)\Vert^2_{L^2(\Omega\setminus\overline{D}}
-\Vert\nabla\epsilon_n\Vert_{L^2(\Omega\setminus\overline{D})}^2
\\
\\
\displaystyle
=\Vert\nabla w\Vert_{L^2(\Omega\setminus\overline{D})}^2
+2\int_{\Omega\setminus\overline{D}}\,\nabla\epsilon_n\cdot\nabla w\,dx
\\
\\
\displaystyle
=-\int_{\partial D}\,\frac{\partial w}{\partial\nu}\,w\,dS+k^2\Vert w\Vert_{L^2(\Omega\setminus\overline{D})}^2
-2\int_{\partial D}\,\epsilon_n\frac{\partial w}{\partial\nu}\,dS
+2\int_{\Omega\setminus\overline{D}}k^2\epsilon w\,dz
\\
\\
\displaystyle
=\int_{\partial D_d}\,\frac{\partial w}{\partial\nu}\,v\,dS
+\int_{\partial D_n}\,\frac{\partial v}{\partial\nu}\,w\,dS
-2\int_{\partial D_d}\,v\frac{\partial w}{\partial\nu}\,dS
\\
\\
\displaystyle
\,\,\,
+k^2\Vert w+\epsilon_n\Vert_{L^2(\Omega\setminus\overline{D})}^2
-k^2\Vert \epsilon_n\Vert_{L^2(\Omega\setminus\overline{D})}^2
\\
\\
\displaystyle
=\left(\int_{\partial D_n}\,\frac{\partial v}{\partial\nu}\,w\,dS
-\int_{\partial D_d}\,\frac{\partial w}{\partial\nu}\,v\,dS\right)
+k^2\Vert w+\epsilon_n\Vert_{L^2(\Omega\setminus\overline{D})}^2
-k^2\Vert \epsilon_n\Vert_{L^2(\Omega\setminus\overline{D})}^2.
\end{array}
\tag {3.15}
$$
Besides we have
$$\begin{array}{l}
\displaystyle
\,\,\,\,\,\,
<(\Lambda_0-\Lambda_D)(v\vert_{\partial\Omega}),v\vert_{\partial\Omega}>
\\
\\
\displaystyle
=-\int_{\partial\Omega}\,\frac{\partial w}{\partial\nu}\,v\,dS\\
\\
\displaystyle
=-\int_{\Omega\setminus\overline{D}}\nabla w\cdot\nabla v\,dz
+\int_{\Omega\setminus\overline{D}}k^2wv\,dz
-\int_{\partial D}\,\frac{\partial w}{\partial\nu}\,v\,dS
\\
\\
\displaystyle
=\int_{\partial D}\,w\frac{\partial v}{\partial\nu}\,dS
-\int_{\partial D}\,\frac{\partial w}{\partial\nu}\,v\,dS
\\
\\
\displaystyle
=-\int_{\partial D_d}\,v\frac{\partial v}{\partial\nu}\,dS
+\int_{\partial D_n}\,w\frac{\partial v}{\partial\nu}\,dS
\\
\\
\displaystyle
\,\,\,
-\int_{\partial D_d}\,\frac{\partial w}{\partial\nu}\,v\,dS
+\int_{\partial D_n}\,\frac{\partial v}{\partial\nu}\,v\,dS\\
\\
\displaystyle
=\left(\int_{\partial D_n}\,\frac{\partial v}{\partial\nu}\,v\,dS
-\int_{\partial D_d}\,v\frac{\partial v}{\partial\nu}\,dS\right)
+\left(\int_{\partial D_n}\,\frac{\partial v}{\partial\nu}\,w\,dS-\int_{\partial D_d}\,\frac{\partial w}{\partial\nu}\,v\,dS
\right).
\end{array}
\tag {3.16}
$$
Here we have
$$\left\{
\begin{array}{l}
\displaystyle
\Vert\nabla v\Vert_{L^2(D_n)}^2-k^2\Vert v\Vert_{L^2(D_n)}^2=\int_{\partial D_n}\,\frac{\partial v}{\partial\nu}\,v\,dS,
\\
\\
\displaystyle
\Vert\nabla v\Vert_{L^2(D_d)}^2-k^2\Vert v\Vert_{L^2(D_d)}^2=\int_{\partial D_d}\,\frac{\partial v}{\partial\nu}\,v\,dS.
\end{array}
\right.
\tag {3.17}
$$
Now from (3.15), (3.16) and (3.17) we obtain (1.16).

Next we have
$$\begin{array}{l}
\displaystyle
\,\,\,\,\,\,
-\Vert\nabla(w+\epsilon_d)\Vert^2_{L^2(\Omega\setminus\overline{D})}
+\Vert\nabla\epsilon_d\Vert_{L^2(\Omega\setminus\overline{D})}^2
\\
\\
\displaystyle
=-\Vert\nabla w\Vert_{L^2(\Omega\setminus\overline{D})}^2
-2\int_{\Omega\setminus\overline{D}}\,\nabla\epsilon_d\cdot\nabla w\,dz
\\
\\
\displaystyle
=\int_{\partial D}\,\frac{\partial w}{\partial\nu}\,w\,dS
+2\int_{\partial D}\,\frac{\partial}{\partial\nu}\epsilon_d\,w\,dS
\\
\\
\displaystyle
\,\,\,
-k^2\Vert w\Vert_{L^2(\Omega\setminus\overline{D})}^2
-2\int_{\Omega\setminus\overline{D}}\,k^2\epsilon_d w\,dz
\\
\\
\displaystyle
=-\int_{\partial D_d}\,\frac{\partial w}{\partial\nu}\,v\,dS
-\int_{\partial D_n}\,\frac{\partial v}{\partial\nu}\,w\,dS
+2\int_{\partial D_n}\,\frac{\partial}{\partial\nu} v\,w\,dS
\\
\\
\displaystyle
\,\,\,
-k^2\Vert w+\epsilon_d\Vert_{L^2(\Omega\setminus\overline{D})}^2
+k^2\Vert\epsilon_d\Vert_{L^2(\Omega\setminus\overline{D})}^2
\\
\\
\displaystyle
=\left(\int_{\partial D_n}\,\frac{\partial v}{\partial\nu}\,w\,dS
-\int_{\partial D_d}\,\frac{\partial w}{\partial\nu}\,v\,dS\right)
-k^2\Vert w+\epsilon_d\Vert_{L^2(\Omega\setminus\overline{D})}^2
+k^2\Vert\epsilon_d\Vert_{L^2(\Omega\setminus\overline{D})}^2.
\end{array}
\tag {3.18}
$$
Thus from (3.16), (3.17) and (3.18) we obtain (1.17).

\noindent
$\Box$

\noindent

Note that once we have found the equation to prove, the proof is just a calculation.
The point of Theorem 3.1 is the introduction of complementing functions $\epsilon_n$ and $\epsilon_d$
in such a way that the integral
$$\displaystyle
\int_{\partial D_n}\frac{\partial v}{\partial\nu}w\,dS
-\int_{\partial D_d}\frac{\partial w}{\partial\nu} v\,dS
$$
has two energy integral expressions given by (3.15) and (3.18).

As a direct corollary, we obtain

\proclaim{\noindent Corollary 3.1.} Let $x\in\Omega$ and
$\sigma\in N_x$.  Let $\{v_m\}$ be an arbitrary needle sequence for $(x,\sigma)$ based on ${\cal G}$.
We have
$$\begin{array}{l}
\displaystyle
\,\,\,\,\,\,
<(\Lambda_0-\Lambda_D)(v_m\vert_{\partial\Omega}),v_m\vert_{\partial\Omega}>
\\
\\
\displaystyle
=\Vert\nabla v_m\Vert_{L^2(D_n)}^2-k^2\Vert v_m\Vert_{L^2(D_n)}^2
\\
\\
\displaystyle
\,\,\,
+\Vert \nabla(w_m+(\epsilon_m)_n)\Vert_{L^2(\Omega\setminus\overline{D})}^2
-k^2\Vert w_m+(\epsilon_m)_n\Vert_{L^2(\Omega\setminus\overline{D})}^2
\\
\\
\displaystyle
\,\,\,
-\Vert \nabla(\epsilon_m)_n\Vert_{L^2(\Omega\setminus\overline{D})}^2
+k^2\Vert (\epsilon_m)_n\Vert_{L^2(\Omega\setminus\overline{D})}^2
\\
\\
\displaystyle
\,\,\,
-\Vert\nabla v_m\Vert_{L^2(D_d)}^2
+k^2\Vert v_m\Vert_{L^2(D_d)}^2
\end{array}
\tag {3.19}
$$
and
$$\begin{array}{ll}
\displaystyle
\,\,\,\,\,\,
<(\Lambda_0-\Lambda_D)(v_m\vert_{\partial\Omega}),v_m\vert_{\partial\Omega}>
\\
\\
\displaystyle
=-\Vert\nabla v_m\Vert_{L^2(D_d)}^2+k^2\Vert v_m\Vert_{L^2(D_d)}^2
\\
\\
\displaystyle
\,\,\,
-\Vert \nabla(w_m+(\epsilon_m)_d)\Vert_{L^2(\Omega\setminus\overline{D})}^2
+k^2\Vert w_m+(\epsilon_m)_d\Vert_{L^2(\Omega\setminus\overline{D})}^2
\\
\\
\displaystyle
\,\,\,
+\Vert \nabla(\epsilon_m)_d\Vert_{L^2(\Omega\setminus\overline{D})}^2
-k^2\Vert (\epsilon_m)_d\Vert_{L^2(\Omega\setminus\overline{D})}^2
\\
\\
\displaystyle
\,\,\,
+\Vert\nabla v_m\Vert_{L^2(D_n)}^2
-k^2\Vert v_m\Vert_{L^2(D_n)}^2,
\end{array}
\tag {3.20}
$$
where $w_m$, $(\epsilon_m)_n$ and $(\epsilon_m)_d$ are given by 
$w_x$, $(\epsilon_x)_n$ and $(\epsilon_x)_d$ with  $G(y,x)$ in (1.6), (1.10) and (1.12) replaced by $v_m(y)$,
respectively.

\endproclaim

\noindent
It should be emphasized that the point is the idea or the principle of the derivation of the things to be proved, like (1.16)
and (1.17) or, (3.19) and (3.20).  It is not a trivial fact as we have already seen.  
It is based on the {\it correspondence principle} mentioned below.

$\quad$

{\bf\noindent Principle.}
Replace the singular solution $G(z,x)$ appeared in some identity, say (3.10) and (3.11), involving
$w_x$, $(\epsilon_x)_n$ and $(\epsilon_x)_d$ with $\{v_m\}$ based on ${\cal G}$.  Then one gets a corresponding
identity for $\{v_m\}$ (to be proved independently), say (3.19) and (3.20).

$\quad$

\noindent
Note that conversely (3.19) and (3.20) yield immediately (3.10) and (3.11) by taking the limit and the formula
(3.1), respectively.

\subsection{Side B of probe method}

The Side B of the probe method is concerned with the blowing up property
of the indicator sequence.  It is based on Theorem 1.2 or Colloary 3.1 and the blowing up property
of the needle sequence stated below.

\proclaim{\noindent Proposition 3.2.}
Given an arbitrary point $x\in\Omega$ and needle $\sigma\in N_x$
let $\{v_m\}$ be an arbitrary needle sequence for $(x,\sigma)$ based on ${\cal G}$.

\noindent 
(a)  Le $V$ be an arbitrary finite cone with vertex at $x$.  Then, we have
$$\displaystyle
\lim_{m\rightarrow\infty}\,\Vert\nabla v_m\Vert_{L^2(V\cap\Omega)}^2=\infty.
$$

\noindent
(b)  Let $z\in\Omega$ be an arbitrary point on $\sigma\setminus\{x\}$ and $B$ open ball centered at $z$.
Then, we have
$$\displaystyle
\lim_{m\rightarrow\infty}\,\Vert\nabla v_m\Vert_{L^2(B\cap\Omega)}^2=\infty.
$$

\endproclaim

\noindent
This fact has been already established in \cite{INew}.

Once we have (3.19), (3.20) and Proposition 3.2, one gets the following theorem which states the
Side B of the probe method.

\proclaim{\noindent Theorem 3.2.}  Let $k=0$.
Let $x\in\Omega$ and $\sigma\in N_x$.
Assume that one of the two cases (a) and (b) listed below are satisfied.

\noindent
(a) $x\in\overline{D}$;

\noindent
(b) $x\in\Omega\setminus\overline{D}$ and $\sigma\cap D\not=\emptyset$.

\noindent
Then for any needle sequence $\{v_m\}$ for $(x,\sigma)$ based on ${\cal G}$ we have
$$\displaystyle
\lim_{m\rightarrow\infty}
<(\Lambda_0-\Lambda_D)(v_m\vert_{\partial\Omega}),v_m\vert_{\partial\Omega}>
=
\left\{
\begin{array}{ll}
\displaystyle
\infty
&
\text{if $\sigma\cap\overline{D_d}=\emptyset$,}
\\
\\
\displaystyle
-\infty
&
\text{if $\sigma\cap\overline{D_n}=\emptyset$.}
\end{array}
\right.
$$
\endproclaim
{\it\noindent Proof.}
We describe only the case when $\sigma\cap\overline{D_d}=\emptyset$.
In this case we have $x\in\Omega\setminus\overline{D_d}$ and the convergence $v_m\rightarrow G(\,\cdot\,,x)$ in $H^2(D_d)$ yields
the boundedness of the sequence $\{(\epsilon_m)_n\}$ in $H^2(\Omega\setminus\overline{D})$.
Thus it follows from (3.19) that
$$\displaystyle
<(\Lambda_0-\Lambda_D)(v_m\vert_{\partial\Omega}),v_m\vert_{\partial\Omega}>
\ge \Vert\nabla v_m\Vert_{L^2(D_n)}-C,
$$
where $C$ is a positive constant.  Here, by Proposition 3.2 under (a) and (b) above we have
$\Vert\nabla v_m\Vert_{L^2(D_n)}\rightarrow\infty$.  This completes the proof.

\noindent
$\Box$

$\quad$

{\bf\noindent Remark 3.3.}  Note that Theorem 3.2 does not cover all the possible cases for $(x,\sigma)$.
For example, if both of the conditions $\sigma\cap\overline{D_d}\not=\emptyset$ and $x\in\overline{D_n}$ are satified,
it would be difficult to state something about the behaviour of the indicator sequence.

$\quad$

The problem is the case when $k\not=0$.  For this, even in the case when $D_d=\emptyset$ we have
only a result  in \cite{INew} under a smalness condition on $k$.  See also Section 2.3.1 of \cite{IReview}
for a concise explanation.

Here, applying the idea described therein to formulae (3.19) and (3.20), we show a result.

We assume that $D_n$ and $D_d$ have the form
$$\begin{array}{ll}
\displaystyle
D_n=\cup_{j=1}^N\,D_{n,\,j}
&
D_d=\cup_{l=1}^M\,D_{d,\,l},
\displaystyle
\end{array}
$$
where $D_{n,\,j}$, $j=1,\cdots,N$ and $D_{d,\,l}$, $l=1,M$ are connected components of $D_n$ and $D_d$, respectively
and satisfiy $\overline{D_{n,\,j}}\cap\overline{D_{n,\,j'}}=\emptyset$ if $j\not=j'$;
 $\overline{D_{d,\,l}}\cap\overline{D_{d,\,l'}}=\emptyset$ if $l\not=l'$.

The assumption on $k$ is as follows: $k$ satisfies all the inequalities listed below:
$$\displaystyle
C(\Omega\setminus\overline{D})^2k^2\le 1;
\tag {3.21}
$$
$$\displaystyle
\max_{j=1,\cdots,N}\,8C(D_{n,j})^2k^2<1;
\tag {3.22}
$$
$$\displaystyle
\max_{l=1,\cdots,M}\,8C(D_{d,j})^2k^2<1.
\tag {3.23}
$$
Here the constants $C(\Omega\setminus\overline{D})$, $C(D_{n,j})$ and $C(D_{d,l})$ denote the Poincar\'e constants \cite{Z}
in the following sense, respectively.

\noindent
(i)  Constant $C(\Omega\setminus\overline{D})$ satisfies, for all $w\in H^1(\Omega\setminus\overline{D})$ with $w=0$ on $\partial\Omega$
$$\displaystyle
\Vert w\Vert_{L^2(\Omega\setminus\overline{D})}
\le C(\Omega\setminus\overline{D})\Vert\nabla w\Vert_{L^2(\Omega\setminus\overline{D})}.
$$

\noindent
(ii)  Constant $C(D_{n,j})$ satisfies, for all $v\in H^1(D_{n,j})$ with $\int_{D_{n,j}}\,vdz=0$
$$\displaystyle
\Vert v\Vert_{L^2(D_{n,j})}
\le C(D_{n,j})\Vert\nabla v\Vert_{L^2(D_{n,j})}.
$$

\noindent
(iii)  Constant $C(D_{d,l})$ satisfies, for all $v\in H^1(D_{d,l})$ with $\int_{D_{d,l}}\,vdz=0$
$$\displaystyle
\Vert v\Vert_{L^2(D_{d,l})}
\le C(D_{d,l})\Vert\nabla v\Vert_{L^2(D_{d,l})}.
$$

\proclaim{\noindent Theorem 3.3.}
All the statements of Theorem 3.2 for $k\not=0$ is valid under the smallness condition (3.21), (3.22) and (3.23).

\endproclaim
{\it\noindent Proof.}
Consider tha case $\sigma$ satisfy $\sigma\cap\overline{D_d}=\emptyset$.
Applying the inequality (3.21) to the function $w=w_m+(\epsilon)_m$, we have
$$\displaystyle
\Vert \nabla(w_m+(\epsilon_m)_n)\Vert_{L^2(\Omega\setminus\overline{D})}^2
-k^2\Vert w_m+(\epsilon_m)_n\Vert_{L^2(\Omega\setminus\overline{D})}^2
\ge 0.
$$
Thus (3.19) yields
$$\begin{array}{l}
\displaystyle
\,\,\,\,\,\,
<(\Lambda_0-\Lambda_D)(v_m\vert_{\partial\Omega}),v_m\vert_{\partial\Omega}>
\\
\\
\displaystyle
\ge \left(\Vert\nabla v_m\Vert_{L^2(D_n)}^2-k^2\Vert v_m\Vert_{L^2(D_n)}^2\right)
+R_m,
\end{array}
\tag {3.24}
$$
where
$$
\displaystyle
R_m=
\displaystyle
-\Vert \nabla(\epsilon_m)_n\Vert_{L^2(\Omega\setminus\overline{D})}^2
-\Vert\nabla v_m\Vert_{L^2(D_d)}^2.
$$
By the convergence property of $\{v_m\}$ in $H^2(D_d)$, we have $\{\epsilon_m\}$ is boundeed in $H^2(\Omega\setminus\overline{D})$.
Therefore the sequence $\{R_m\}$ is bounded.  
Besides, applying the same argument in \cite{INew} (and also see \cite{IReview}) to the first term of  the right-hand side on 
(3.24), we obtain
$$\displaystyle
\Vert\nabla v_m\Vert_{L^2(D_n)}^2-k^2\Vert v_m\Vert_{L^2(D_n)}^2
\ge C_1\Vert\nabla v_m\Vert_{L^2(D_n)}^2 -C_2,
$$
here $C_1$ and $C_2$ are positive constants independent of $m$, however, depends on $\sigma$, $D_n$ and $k$ satisfying (3.22).
Thus the blowing up property of the indicator sequence is reduced to that of $\Vert\nabla v_m\Vert_{L^2(D_n)}^2$, that is covered by Proposition 3.2.

The treatement of the case when $\sigma\cap\overline{D_n}=\emptyset$ is the same except for the use of (3.23) and  (3.20) instead of  (3.22) and 
(3.19), respectively.
\noindent
$\Box$

$\quad$

{\bf\noindent Remark 3.4.}
It should be noted that, in \cite{CLNW}  they considered the probe method \cite{Iwave,I2} for the Helmholtz equation
$\Delta u+k^2 u=0$
in the mixed obstacle case.
However, in their paper only the Side A of the probe method is considered
and their argument is based on a combination of that of \cite{Iwave} and a detailed singularity analysis of the reflected solution.
There is is no description about the Side B of the probe method, which has been 
introduced in \cite{INew} and developed in \cite{IHokkai}.
Besides,  even the case when the wave number $k=0$ their
result does not cover Theorem 3.2.  This is due to the lack of formulae (1.16) and (1.17)
or (3.19) and (3.20).

\subsection{Singular sources method included in IPS}

The singular sources method consists of three parts listed below.

$\quad$

\noindent
(a)  Given $x\in\Omega\setminus\overline{D}$ and $\sigma\in N_x$
let $\{v_n\}$ be an arbitrary needle sequence for $(x,\sigma)$ based on ${\cal G}$.
Then we have formula (3.6), that is,
$$
\displaystyle
w_x(x)
=-\lim_{n\rightarrow\infty}
<(\Lambda_0-\Lambda_D)(v_n\vert_{\partial\Omega}),
(G(\,\cdot\,,x)-v_n)\vert_{\partial\Omega}>.
\tag {3.25}
$$

$\quad$

\noindent
(b)  It holds that:

\noindent
(i)  $\lim_{x\rightarrow a\in\partial D_n}w_x(x)=\infty$;

\noindent
(ii)  $\lim_{x\rightarrow b\in\partial D_d}w_x(x)=-\infty$.

$\quad$

\noindent
(c)  For each $\epsilon_i>0$, $i=1,2$
$$\displaystyle
\sup_{x\in\Omega\setminus\overline{D},\,\text{dist}\,(x,\,\partial D)>\epsilon_1,
\,\text{dist}\,(x,\,\partial\Omega)>\epsilon_2}\,\vert w_x(x)\vert<\infty.
\tag {3.26}
$$

$\quad$

\noindent
The statements (b) and (c) are the direct consequence of  Corollary 1.1, outer decomposition (1.8),  (1.3) and (1.7).

Note that Remark 3.1 works also for (3.26) in the case when ${\cal G}={\cal G}^0$.  That is,
the condition $\text{dist}\,(x,\partial\Omega)>\epsilon_2$ in (3.26)  is dropped.

$\quad$

{\bf\noindent Remark 3.5.}
If ${\cal G}={\cal G}^*$, then $G(\,\cdot\,,x)=0$ on $\partial\Omega$ and  by (3.25) one gets
$$\displaystyle
w_x(x)
=\lim_{n\rightarrow\infty}
<(\Lambda_0-\Lambda_D)(v_n\vert_{\partial\Omega}), v_n\vert_{\partial\Omega}>
$$
provided $\sigma\cap\overline{D}=\emptyset$.  From this together with  (3.1) and (3.2)
we obtain
$$\displaystyle
w_x(x)=I(x)=W_x(x).
$$
So this is the {\it completely integrated version} of the probe and singular sources methods.
By Theorem 3.2 this version also has the Side B.
To distinguish from other cases we denote $w_x(x)=I(x)=W_x(x)$
by $w_x^*(x)=I^*(x)=W_x^*(x)$ if ${\cal G}={\cal G}^*$.

\subsection{Side B of singular sources method}

Given $x\in\Omega$ and $\sigma\in N_x$ 
let $\{v_n^0\}$ be the needle sequence for $(x,\sigma)$ based on ${\cal G}={\cal G}^0$, that is
$$\displaystyle
v_n^0\rightarrow G(\,\cdot\,-x)
$$
in $H^2_{\text{loc}}(\Omega\setminus\sigma)$.
Let $H(z)=H(z,x)$ solve 
$$
\left\{
\begin{array}{ll}
\displaystyle
\Delta H+k^2H=0, & z\in\Omega,\\
\\
\displaystyle
H(z)=-G(z-x), & z\in\partial\Omega.
\end{array}
\right.
$$
It is clear that the $H(\,\cdot\,,x)$ satisfies (1.3).
The function 
$$\begin{array}{ll}
\displaystyle
v_n(z)=v_n^0(z)+H(z,x), & z\in\Omega,
\end{array}
$$
satisfies the Helmholtz equation in $\Omega$
and that the sequence $\{v_n\}$ satisfies
$$\displaystyle
v_n\rightarrow G(\,\cdot\,-x)+H(\,\cdot\,,x)
$$
in $H^2_{\text{loc}}(\Omega\setminus\sigma)$.  This means that sequence 
$\{v_n\}$ is a needle sequence for $(x,\sigma)$ based on ${\cal G}={\cal G}^*$
(see also Remark 1.1).
Thus, if $\sigma\cap\overline{D}=\emptyset$, then
by Remark 3.5 we have
$$\displaystyle
w_x^*(x)=I^*(x)=W_x^*(x)
=\lim_{n\rightarrow\infty}
<(\Lambda_0-\Lambda_D)(v_n\vert_{\partial\Omega}), v_n\vert_{\partial\Omega}>.
$$
Here note that we have
$$\begin{array}{ll}
\displaystyle
v_n(z)=v_n^0(z)-G(z-x), & z\in\partial\Omega.
\end{array}
$$
Therefore we obtain
$$\displaystyle
w_x^*(x)=I^*(x)=W_x^*(x)=\lim_{n\rightarrow\infty}
<(\Lambda_0-\Lambda_D)((v_n^0-G(\,\cdot\,-x))\vert_{\partial\Omega}), 
(v_n^0-G(\,\cdot\,-x))\vert_{\partial\Omega}>.
$$
Besides, as a corollary of Theorem 3.2 we obtain

\proclaim{\noindent Corollary 3.2.}  Let $k=0$.
Let $x\in\Omega$ and $\sigma\in N_x$.
Assume that one of the two cases (a) and (b) listed in Theorem 3.2 are satisfied.
\noindent
Then for any needle sequence $\{v_m^0\}$ for $(x,\sigma)$ based on ${\cal G}={\cal G}^0$ we have
$$\displaystyle
\lim_{n\rightarrow\infty}
<(\Lambda_0-\Lambda_D)((v_n^0-G(\,\cdot\,-x))\vert_{\partial\Omega}),
(v_n^0-G(\,\cdot\,-x))\vert_{\partial\Omega}>
=
\left\{
\begin{array}{ll}
\displaystyle
\infty
&
\text{if $\sigma\cap\overline{D_d}=\emptyset$,}
\\
\\
\displaystyle
-\infty
&
\text{if $\sigma\cap\overline{D_n}=\emptyset$.}
\end{array}
\right.
$$
\endproclaim

And also as a corollary of Theorem 3.3 we have

\proclaim{\noindent Corollary 3.3.}
Let $k\ge 0$ satisfy (3.21), (3.22) and (3.23).  Then the same conclusions as Corollary 3.2 are valid.

\endproclaim

$\quad$

Let $w_x^0(y)=w_x(y;{\cal G}^0)$.
The $w_x^0$ solves
$$\left\{
\begin{array}{ll}
\displaystyle
\Delta w+k^2w=0, & y\in\Omega\setminus\overline{D},\\
\\
\displaystyle
\frac{\partial w}{\partial\nu}=-\frac{\partial}{\partial\nu}G(y-x), & y\in\partial D_n,\\
\\
\displaystyle
w=-G(y-x), & y\in\partial D_d,\\
\\
\displaystyle
w=0,  & y\in\partial\Omega.
\end{array}
\right.
$$
The function $\Omega\setminus\overline{D}\ni x\mapsto w_x^0(x)$
is a natural extension of  the {\it indicator function for the singular sources method} discussed therein to the case $D_d\not=\emptyset$.
See also Subsection 1.3 in \cite{IPS}
for an explanation of why $w_x^0(x)$ is called the indicator function for the singular sources method
in relation to its original singular sources method of Potthast \cite{P1}.

Now let $\sigma\cap\overline{D}=\emptyset$.
By (3.4) and (3.25) in the case ${\cal G}={\cal G}^0$, we have
$$
\displaystyle
w_x^0(x)
=-\lim_{n\rightarrow\infty}<(\Lambda_0-\Lambda_D)v_n^0\vert_{\partial\Omega},G_n(\,\cdot\,,x)\vert_{\partial\Omega}>,
$$
where
$$\displaystyle
G_n(z,x)=G(z-x)-v_n^0(z).
$$
Here we have the trivial decomposition
$$\begin{array}{l}
\displaystyle
\,\,\,\,\,\,
<(\Lambda_0-\Lambda_D)G(\,\cdot\,-x)\vert_{\partial\Omega},G(\,\cdot\,-x)\vert_{\partial\Omega}>
\\
\\
\displaystyle
=<(\Lambda_0-\Lambda_D)v_n^0\vert_{\partial\Omega},v_n^0\vert_{\partial\Omega}>
+<(\Lambda_0-\Lambda_D)G_n(\,\cdot\,,x)\vert_{\partial\Omega},G_n(\,\cdot\,,x)\vert_{\partial\Omega}>
\\
\\
\displaystyle
\,\,\,
+<(\Lambda_0-\Lambda_D)v_n^0\vert_{\partial\Omega},G_n(\,\cdot\,,x)\vert_{\partial\Omega}>
+<(\Lambda_0-\Lambda_D)(G_n(\,\cdot\,,x)\vert_{\partial\Omega}), v_n^0\vert_{\partial\Omega}>.
\end{array}
$$
This together with the symmetry of Dirichlet-to-Neumann maps $\Lambda_0$ and $\Lambda_D$
yields the expression
$$\begin{array}{l}
\displaystyle
\,\,\,\,\,\,
-<(\Lambda_0-\Lambda_D)v_n^0\vert_{\partial\Omega},G_n(\,\cdot\,,x)\vert_{\partial\Omega}>
\\
\\
\displaystyle
=\frac{1}{2}
\left(<(\Lambda_0-\Lambda_D)v_n^0\vert_{\partial\Omega},v_n^0\vert_{\partial\Omega}>
+<(\Lambda_0-\Lambda_D)G_n(\,\cdot\,,x)\vert_{\partial\Omega},G_n(\,\cdot\,,x)\vert_{\partial\Omega}>\right)\\
\\
\displaystyle
\,\,\,
-\frac{1}{2}<(\Lambda_0-\Lambda_D)G(\,\cdot\,-x)\vert_{\partial\Omega},G(\,\cdot\,-x)\vert_{\partial\Omega}>.
\end{array}
\tag {3.27}
$$
Therefore, using Theorem 3.2  for the choice ${\cal G}={\cal G}^0$ and Corollary 3.2, we obtain the side $B$ of the singular sources method formulated in \cite{IPS}.

\proclaim{\noindent Corollary 3.4.}  Let $k=0$.
Let $x\in\Omega$ and $\sigma\in N_x$.
Assume that one of the two cases (a) and (b) listed in Theorem 3.2 are satisfied.
\noindent
Then for any needle sequence $\{v_n^0\}$ for $(x,\sigma)$ based on ${\cal G}={\cal G}^0$ we have
$$\displaystyle
-\lim_{n\rightarrow\infty}
<(\Lambda_0-\Lambda_D)v_n^0\vert_{\partial\Omega},(G(\,\cdot\,-x)-v_n^0)\vert_{\partial\Omega}>
=
\left\{
\begin{array}{ll}
\displaystyle
\infty
&
\text{if $\sigma\cap\overline{D_d}=\emptyset$,}
\\
\\
\displaystyle
-\infty
&
\text{if $\sigma\cap\overline{D_n}=\emptyset$.}
\end{array}
\right.
$$
\endproclaim

And also from Theorem 3.3 and Corollary 3.3 we have

\proclaim{\noindent Corollary 3.5.}
Let $k\ge 0$ satisfy (3.21), (3.22) and (3.23).  Then we have the same conclusions as Corollary 3.4.

\endproclaim

\noindent
Corollaries 3.4 and 3.5 could never be obtained using a single methodology, and show us the greatest advantage of the integrated theory.

$\quad$

{\bf\noindent Remark 3.6.}
It follows from (3.27) and the existence of the needle sequence
which is a consequence of the Runge approximation property for the Helmholtz equation in $\Omega$ we have the expression
$$\displaystyle
w_x^0(x)=\frac{1}{2}(I^0(x)+I^*(x)-<(\Lambda_0-\Lambda_D)G(\,\cdot\,-x)\vert_{\partial\Omega},G(\,\cdot\,-x)\vert_{\partial\Omega}>),
$$
where the $I^0(x)$ denotes the $I(x)$ given by (3.2) (or both of  (3.10) and  (3.11)) with the case when ${\cal G}={\cal G}^0$, that is, $G(\,\cdot\,,x)=G(\,\cdot\,-x)$.

\subsection{Additional remarks}

\subsubsection{Lifting}

First for general ${\cal G}$, from Lemma 2.1 we obtain
$$\displaystyle
w_x(y)-\int_{\partial\Omega}\frac{\partial}{\partial\nu}w_x(z)G(z,y)\,dS(z)
=w_y(x)-\int_{\partial\Omega}\frac{\partial}{\partial\nu}w_y(z)G(z,x)\,dS(z).
\tag {3.28}
$$
So define
$$\begin{array}{ll}
\displaystyle
I(x,y)=w_x(y)-\int_{\partial\Omega}\frac{\partial}{\partial\nu}w_x(z)G(z,y)\,dS(z), & (x,y)\in(\Omega\setminus\overline{D})^2.
\end{array}
\tag {3.29}
$$
Then (3.28) yields the symmetry
$$\displaystyle
I(x,y)=I(y,x).
$$
Besides, by (3.1) and (3.3), we have another expression for the idicator function.
$$\displaystyle
I(x)=w_x(x)-\int_{\partial\Omega}\frac{\partial}{\partial\nu}w_x(z)G(z,x)\,dS(z).
$$
This together with (3.29) yields $I(x)=I(x,y)\vert_{y=x}$ and in this sense, the $I(x,y)$ is called the {\it lifting} of $I(x)$.

The inner decomposition (3.14) itself has the lifted version.
For general ${\cal G}$, by  (2.12) we have
$$\begin{array}{ll}
\displaystyle
w_x^1(y)=I^1(x,y)-\int_{\partial\Omega}\,G(z,x)\frac{\partial}{\partial\nu}w_y(z)\,dS(z),
&
\displaystyle
(x,y)\in (\Omega\setminus\overline{D})^2,
\end{array}
\tag {3.30}
$$
where
$$\displaystyle
I^1(x,y)=\int_{\Omega\setminus\overline{D}}\nabla w_x^1(z)\cdot\nabla w_y^1(z)\,dz-\int_{\Omega\setminus\overline{D}}k^2 w_x^1(z)w_y^1(z)\,dz
-\int_{\partial\Omega}\,\frac{\partial}{\partial\nu}G(z,y)G(z,x)\,dS(z).
$$
A similar computation to (3.8) yields
$$\displaystyle
<\Lambda_D(G(\,\cdot\,,y)\vert_{\partial\Omega}),G(\,\cdot\,,x)\vert_{\partial\Omega}>
=\int_{\Omega\setminus\overline{D}}\nabla w_y^1(z)\cdot\nabla w_x^1(z)\,dz-\int_{\Omega\setminus\overline{D}}k^2 w_y^1(z)w_x^1(z)\,dz.
$$
Thus we have
$$\displaystyle
I^1(x,y)=<\Lambda_D(G(\,\cdot\,,y)\vert_{\partial\Omega}),G(\,\cdot\,,x)\vert_{\partial\Omega}>
-\int_{\partial\Omega}\,\frac{\partial}{\partial\nu}G(z,y)G(z,x)\,dS(z).
\tag {3.31}
$$
By (3.13), this yields $I^1(x)=I^1(x,y)\vert_{y=x}$.  Thus, $I^1(x,y)$ gives  a lifting of $I^1(x)$.

Next rewrite (3.29) as
$$\displaystyle
w_x(y)=I(x,y)+\int_{\partial\Omega}G(z,y)\frac{\partial}{\partial\nu}w_x(z)\,dS(z).
$$
This together with (3.30) yields
$$\begin{array}{ll}
\displaystyle
W_x(y)
&
\displaystyle
=I(x,y)+I^1(x,y)
\\
\\
\displaystyle
&
\displaystyle
\,\,\,
+\int_{\partial\Omega}G(z,y)\frac{\partial}{\partial\nu}w_x(z)\,dS(z)-\int_{\partial\Omega}\,G(z,x)\frac{\partial}{\partial\nu}w_y(z)\,dS(z).
\end{array}
\tag {3.32}
$$
This is the lifted version of inner decomposition (3.14).  
Note also that we have {\it twisted} decomposition
$$\displaystyle
w_y(x)+w_x^1(y)=I(x,y)+I^1(x,y).
\tag {3.33}
$$

\subsubsection{Uniqueness}

We consider only two cases: ${\cal G}={\cal G}^0, {\cal G}^*$.
From (3.29) we have
$$\displaystyle
I(x,y)=
\left\{
\begin{array}{ll}
\displaystyle
w_x(y)-\int_{\partial\Omega}\frac{\partial}{\partial\nu}w_x(z)G(z-y)\,dS(z), & \text{if ${\cal G}={\cal G}^0$},
\\
\\
\displaystyle
w_x(y), & \text{if ${\cal G}={\cal G}^*$.}
\end{array}
\right.
\tag {3.34}
$$
Note that $w_x(y)$ depends on ${\cal G}$ and the symmetry of $I(x,y)$ yields the symmetry $w_x(y)=w_y(x)$ in the case ${\cal G}={\cal G}^*$.

The expression (3.34) together with symmetry of $I(x,y)$ yields:

\noindent
for each fixed $x\in\Omega\setminus\overline{D}$ 
$$\begin{array}{ll}
\displaystyle
\Delta_yI(x,y)+k^2 I(x,y)=0, & y\in\Omega\setminus\overline{D};
\end{array}
$$

\noindent
for each fixed $y\in\Omega\setminus\overline{D}$ 
$$\begin{array}{ll}
\displaystyle
\Delta_xI(x,y)+k^2 I(x,y)=0, & x\in\Omega\setminus\overline{D}.
\end{array}
$$
Here the symbols $\Delta_y$ and $\Delta_x$ denote the Laplacian with respect to $y$ and $x$, respectively.

By the unique continuation property of the Helmholtz equation, one concludes:
indicator function $I(x)$, $x\in\Omega\setminus\overline{D}$ is uniquely determined by $I(x,y)$ given at all $(x,y)\in U\times V$,
where $U$ and $V$ are arbitrary nonempty open subsets of $\Omega\setminus\overline{D}$, typically in a small neighbourhood of $\partial\Omega$.
Besides, using the argument for the proof of (3.3), we obtain also the computation formula of the lifting
$$\displaystyle
I(x,y)=\lim_{n\rightarrow\infty}<(\Lambda_0-\Lambda_D)(v_n\vert_{\partial\Omega}),v_n'\vert_{\partial\Omega}>,
$$
where $v_n$ is the same as that of (3.3) and $v_n'$ is an arbitrary needle sequence for $(y,\sigma')$ based on ${\cal G}$ and $\sigma'\in N_y$
satisfying $\sigma'\cap\overline{D}=\emptyset$.  So {\it in principle} or theoretically, it suffices to use only the needle sequences for 
the needles with tips in $U\times V$, say with $U=V$ and $U$ is given by the intersection of  a {\it small} open ball centered at a point on $\partial\Omega$
with $\Omega$.  In that case we can use only the straight needles {\it explicitly constructed} in \cite{ICar}.

Summing up, we have obtained the following {\it uniqueness theorem} by using needles localized, say in a {\it small} neighbourhood of $\partial\Omega$ in $\overline{\Omega}$.

\proclaim{\noindent Proposition 3.1.}   Let ${\cal G}={\cal G}^0, {\cal G}^*$.  
Let $U$ be an arbitray nonempty open subset of $\Omega\setminus\overline{D}$.   Assume that we have the data $\Lambda_D(v_n\vert_{\partial\Omega})$
for all $x\in U$ and a needle $\sigma\in N_x$ with $\sigma\setminus\partial\Omega\subset U$,
and  a needle sequence $\{v_n\}$ for $(x,\sigma)$ based on ${\cal G}$.
Then, the obstacles
$D_d$ and $D_n$ are uniquely determined by the data.
\endproclaim

\noindent
The key of the proof  is to put the calculation process of the lifting $I(x,y)$ for all $(x, y)\in U^2$ in between.
This result could never have been found using a single methodology alone.

\subsubsection{Symmetry of $I^1(x,y)$ and implications}

For general ${\cal G}$ from (3.32) we obtain
$$\begin{array}{ll}
\displaystyle
\frac{W_x(y)+W_y(x)}{2}
&
\displaystyle
=I(x,y)+\frac{I^1(x,y)+I^1(y,x)}{2}.
\end{array}
\tag {3.35}
$$
Note that for general ${\cal G}$ the $I^1(x,y)$ is not necessary symmetric with respect to variables $x$ and $y$.  
Here we note that the $I^1(x,y)$ is symmetric for ${\cal G}={\cal G}^*, {\cal G}^0$. 
In fact, if ${\cal G}={\cal G^*}$, then $G(\,\cdot\,,x)=0$ on $\partial\Omega$ for each $x\in\Omega\setminus\overline{D}$.
Thus (3.31) yields $I^1(x,y)=0=I^1(y,x)$.

For general ${\cal G}$, a similar argument for the proof of the symmetry of Green's function, we have, for $(x,y)\in\Omega^2$ with $x\not=y$
$$\begin{array}{l}
\displaystyle
\,\,\,\,\,\,
\int_{\partial\Omega}\,\frac{\partial}{\partial\nu}G(z,y)G(z,x)\,dS(z)
\\
\\
\displaystyle
=-G(y,x)+\int_{\Omega}\nabla G(z,y)\cdot\nabla G(z,x)\,dz
-\int_{\Omega}k^2 G(z,y)G(z,x)\,dz.
\end{array}
$$
Note that all the integrands are absolutley integrable since $x\not=y$.    Thus $I^1(x,y)$ with $x\not=y$ takes the form
$$\begin{array}{ll}
\displaystyle
I^1(x,y)
&
\displaystyle
=\int_{\Omega\setminus\overline{D}}\nabla w_x^1(z)\cdot\nabla w_y^1(z)\,dz-\int_{\Omega\setminus\overline{D}}k^2 w_x^1(z)w_y^1(z)\,dz
\\
\\
\displaystyle
&
\displaystyle
\,\,\,
+G(y,x)-\int_{\Omega}\nabla G(z,y)\cdot\nabla G(z,x)\,dz+\int_{\Omega}k^2 G(z,y)G(z,x)\,dz.
\end{array}
$$
Recalling the expression (1.2), we see that $I^1(x,y)=I^1(y,x)$ if and only if $H(y,x)=H(x,y)$.  Thus $I^1(x,y)$ is symmetric in the case ${\cal G}={\cal G}^0$ since $H\equiv 0$.

Therefore from (3.35) we obtain, for ${\cal G}={\cal G}^0, {\cal G}^*$\footnote{Note that, in particular, if ${\cal G}={\cal G}^*$ we have $W_x(y)=I(x,y)=w_x(y)$. This implies the symmetry of $W_x(y)=w_x(y)$.}
$$\displaystyle
\frac{W_x(y)+W_y(x)}{2}=I(x,y)+I^1(x,y).
\tag {3.36}
$$
And this together with (3.33) yields {\it twisted symmetry}:
$$\displaystyle
w_x(y)-w_y^1(x)=w_y(x)-w_x^1(y).
$$
Besides, the expression (3.31)  together with symmetry yields:

\noindent
for each fixed $y\in\Omega\setminus\overline{D}$ we have
$$\begin{array}{ll}
\displaystyle
\Delta_x I^1(x,y)+k^2 I^1(x,y)=0, & x\in\Omega\setminus\overline{D};
\end{array}
$$

\noindent
for each fixed $x\in\Omega\setminus\overline{D}$, we have
$$\begin{array}{ll}
\displaystyle
\Delta_y I^1(x,y)+k^2 I^1(x,y)=0, & y\in\Omega\setminus\overline{D}.
\end{array}
$$

\noindent
Besides, using (3.33) we conclude that, for each fixed $y\in\Omega\setminus\overline{D}$
$$\begin{array}{ll}
\displaystyle
\Delta_x(w_x^1(y))+k^2w_x^1(y)=0, & x\in\Omega\setminus\overline{D}
\end{array}
$$
and for each fixed $x\in\Omega\setminus\overline{D}$
$$\begin{array}{ll}
\displaystyle
\Delta_y (w_y(x))+k^2 w_y(x)=0, & y\in\Omega\setminus\overline{D}.
\end{array}
$$
Finally from (3.36) we obtain, for each fixed $x\in\Omega\setminus\overline{D}$
$$\begin{array}{ll}
\displaystyle
\Delta_y (W_y(x))+k^2 W_y(x)=0, & y\in\Omega\setminus\overline{D}.
\end{array}
$$
As a conclusion, we have:
\proclaim{\noindent Proposition  3.2.}
Let ${\cal G}={\cal G}^0, {\cal G}^*$.  
Let $U$ be an arbitray nonempty open subset of $\Omega\setminus\overline{D}$.  
Then the values $W_x(x)$ and $w_x(x)$ at all $x\in\Omega\setminus\overline{D}$ are uniquely determined
by those of $W_x(y)$ and $w_x(y)$ for all $(x,y)\in U^2$, respectively.
\endproclaim

\section{Conclusion and remarks}

Now it became clear that the IPS function plays the central role
in deriving the probe and singular sources methods.  Besides, the {\it method of complemening function}, which is introduced in the proof of Theorem 1.1, 
makes everything so clear.   
Everything about the both methods
can be derived from the knowledge of the IPS function $W_x(x)$.
As a byproduct, we found the Side B of both the probe and singular sources methods
for the mixed obstacle case.  This is an advantage of IPS.  However, there is a proviso that this comes at the expense of the regularity of the boundaries of the obstacles and whole domain.  This seems to be unavoidable in order to establish especially the singular sources method since its is based on Green's theorem.

In this paper, we have considered only the case when the governing equation is given by the Helmholtz equation.
However, the method developed here can be applied also to the same type of inverse obstacle problems governed by various partial differential equations, for example,
the Navier equation, the Stokes system, the biharmonic equation, and so on.  
And also it would be interested to consider their time domain versions by the spirit of IPS.  Those belong to our next
project.

Our theory yields also an {\it alternative simple proof} of a result on the probe method described in \cite{CLNW}, which is
nothing but the Side A called in this paper.
However, it should be pointed out that 
the Side B {\it without smallness} of $k$ is still open at the present time even for the case $D_d=\emptyset$, see \cite{IReview}.

$$\quad$$

\centerline{{\bf Acknowledgment}}

The author was supported by Grant-in-Aid for
Scientific Research (C)(No. 24K06812) of Japan  Society for
the Promotion of Science.

$$\quad$$

\section{Appendix}

Note that in this appendix it is assumed that $k^2$ satisfies Assumption 1.

\proclaim{\noindent Lemma A.}  
Let $v\in H^2(\Omega)$ be an arbitrary solution of the Helmholtz equation $\Delta v+k^2v=0$ in $\Omega$.
Let $u\in H^2(\Omega\setminus\overline{D})$ be the solution of  (1.1) with $f=v$ on $\partial\Omega$.
We have
$$\displaystyle
\Vert u-v\Vert_{L^2(\Omega\setminus\overline{D})}\le C\,\left(\Vert v\Vert_{L^2(D_n)}+\Vert v\Vert_{L^2(\partial D_d)}\right),
$$
where $C$ is positive constant independent of $v$.

\endproclaim

{\it\noindent Proof.}
Set $w=u-v$.  The $w=w(y)$ satisfies
$$\left\{
\begin{array}{ll}
\displaystyle
\Delta w+k^2w=0, & y\in\Omega\setminus\overline{D},\\
\\
\displaystyle
\frac{\partial w}{\partial\nu}=-\frac{\partial}{\partial\nu}v(y), & y\in\partial D_n,\\
\\
\displaystyle
w=-v(y), & y\in\partial D_d,\\
\\
\displaystyle
w=0, & y\in\partial\Omega.
\end{array}
\right.
$$
Decompose $w$ as $w=w_1+w_2$, where the $w_1=w_1(y)$ solves
$$\left\{
\begin{array}{ll}
\displaystyle
\Delta w_1+k^2w_1=0, & y\in\Omega\setminus\overline{D},\\
\\
\displaystyle
\frac{\partial w_1}{\partial\nu}=-\frac{\partial}{\partial\nu}v(y), & y\in\partial D_n,\\
\\
\displaystyle
w_1=0, & y\in\partial D_d,\\
\\
\displaystyle
w_1=0, & y\in\partial\Omega
\end{array}
\right.
$$
and thus $w_2=w_2(y)$ satisfies
$$\left\{
\begin{array}{ll}
\displaystyle
\Delta w_2+k^2w_2=0, & y\in\Omega\setminus\overline{D},\\
\\
\displaystyle
\frac{\partial w_2}{\partial\nu}=0, & y\in\partial D_n,\\
\\
\displaystyle
w_2=-v(y), & y\in\partial D_d,\\
\\
\displaystyle
w_2=0, & y\in\partial\Omega.
\end{array}
\right.
$$
Considering $\Omega\setminus\overline{D}$ as $(\Omega\setminus\overline{D_d})\setminus\overline{D_n}$
and applying Lemma 2.2 of \cite{IR} to the case when $\Omega$ and $D$ are replaced with $\Omega\setminus\overline{D_d}$ and $D_n$, respectively, we have
$$\displaystyle
\Vert w_1\Vert_{L^2(\Omega\setminus\overline{D})}\le C\Vert v\Vert_{L^2(D_n)}.
$$
So the problem is to show that
$$\displaystyle
\Vert w_2\Vert_{L^2(\Omega\setminus\overline{D})}\le C\Vert v\Vert_{L^2(\partial D_d)}.
\tag {A.1}
$$
Here we employ a slightly modified argument for the proof of (4.12) in Lemma 4.1 of \cite{IE00}.
Solve
$$\left\{
\begin{array}{ll}
\displaystyle
\Delta p+k^2p=w_2, & y\in\Omega\setminus\overline{D},\\
\\
\displaystyle
\frac{\partial p}{\partial\nu}=0, & y\in\partial D_n,\\
\\
\displaystyle
p=0, & y\in\partial D_d,\\
\\
\displaystyle
p=0, & y\in\partial\Omega.
\end{array}
\right.
$$
Then we have
$$\begin{array}{ll}
\displaystyle
\int_{\Omega\setminus\overline{D}}w_2^2\,dy
&
\displaystyle
=\int_{\Omega\setminus\overline{D}}pw_2\,dy
\\
\\
\displaystyle
&
\displaystyle
=\int_{\Omega\setminus\overline{D}}(\Delta p+k^2p)w_2\,dy
\\
\\
\displaystyle
&
\displaystyle
=-\int_{D_n}\frac{\partial p}{\partial\nu}w_2\,dS
-\int_{D_d}\frac{\partial p}{\partial\nu}w_2\,dS
-\int_{\Omega\setminus\overline{D}}\nabla p\cdot\nabla w_2\,dy+\int_{\Omega\setminus\overline{D}}k^2pw_2dy
\\
\\
\displaystyle
&
\displaystyle
=\int_{\partial D_d}\frac{\partial p}{\partial\nu}v\,dS
+\int_{\partial D_n}\frac{\partial w_2}{\partial\nu}p\,dS
+\int_{\partial D_d}\frac{\partial w_2}{\partial\nu}p\,dS
\\
\\
\displaystyle
&
\displaystyle
=\int_{\partial D_d}\frac{\partial p}{\partial\nu}v\,dS
+\int_{\partial D_d}\frac{\partial w_2}{\partial\nu}p\,dS
\\
\\
\displaystyle
&
\displaystyle
=\int_{\partial D_d}\frac{\partial p}{\partial\nu}v\,dS.
\end{array}
$$
Thus one gets
$$\displaystyle
\Vert w_2\Vert_{L^2(\Omega\setminus\overline{D})}^2
\le 
\Vert \nabla p\Vert_{L^2(\partial D_d)}\Vert v\Vert_{L^2(\partial D_d)}
$$
By elliptic regularity up to boundary, we have
$$\displaystyle
\Vert p\Vert_{H^2(\Omega\setminus\overline{D})}\le C\Vert w_2\Vert_{L^2(\Omega\setminus\overline{D})}
$$
and the trace theorem yields
$$\displaystyle
\Vert \nabla p\Vert_{H^{\frac{1}{2}}(\partial D_d)}
\le C\Vert \nabla p\Vert_{H^1(\Omega\setminus\overline{D})}.
$$
Combining these, we obtain (A.1).

\noindent
$\Box$

\proclaim{\noindent Lemma B.}
The solution $w=w_x$ of (1.6) for $x\in\Omega\setminus\overline{D}$ satisfies, for each $\epsilon>0$
$$\displaystyle
\sup_{x\in\Omega\setminus\overline{D},\,\epsilon<\text{dist}\,(x,\partial\Omega)}\Vert w_x\Vert_{L^2(\Omega\setminus\overline{D})}<\infty.
$$
\endproclaim
{\it\noindent Proof.}
There are two ways to validate the statement.
The first one is a combination of Lemma A and a limiting argument based on the Runge approximation propoerty for the Helmholtz equation
in the whole domain $\Omega$ provided $k^2$ is not a Dirichlet eigenvalue for the minus Laplacian $-\Delta$ in $\Omega$.
Another one goes back to the idea for establishing the estimate (19) and  (28) in \cite{Iwave} for the case when $D_n=\emptyset$ and $D_d=\emptyset$,
respectively.
Since the later one is elementary, we describe here.
Decompose $w_x$ as $w_x=(w_1)_x+(w_2)_x$, where the $(w_1)_x=w_1(y)$ solves
$$\left\{
\begin{array}{ll}
\displaystyle
\Delta w_1+k^2w_1=0, & y\in\Omega\setminus\overline{D},\\
\\
\displaystyle
\frac{\partial w_1}{\partial\nu}=-\frac{\partial}{\partial\nu}G(y,x), & y\in\partial D_n,\\
\\
\displaystyle
w_1=0, & y\in\partial D_d,\\
\\
\displaystyle
w_1=0, & y\in\partial\Omega
\end{array}
\right.
$$
and thus $(w_2)_x=w_2(y)$ satisfies
$$\left\{
\begin{array}{ll}
\displaystyle
\Delta w_2+k^2w_2=0, & y\in\Omega\setminus\overline{D},\\
\\
\displaystyle
\frac{\partial w_2}{\partial\nu}=0, & y\in\partial D_n,\\
\\
\displaystyle
w_2=-G(y,x), & y\in\partial D_d,\\
\\
\displaystyle
w_2=0, & y\in\partial\Omega.
\end{array}
\right.
$$
Considering $\Omega\setminus\overline{D}$ as $(\Omega\setminus\overline{D_d})\setminus\overline{D_n}$
and applying the argument for the proof of (28) in \cite{Iwave} to the case when $\Omega$ and $D$ are replaced with $\Omega\setminus\overline{D_d}$ and $D_n$, respectively, we have
$$\displaystyle
\Vert(w_1)_x\Vert_{L^2(\Omega\setminus\overline{D})}
\le
C\left(\int_{\partial D_n}\vert z-x\vert^{\frac{1}{2}}\,\left\vert\frac{\partial G}{\partial\nu}(z,x)\right\vert\,dz
+\int_{D_n}\vert G(z,x)\vert\,dz\right),
$$
where $C$ is a positive constant independent of $x\in\Omega\setminus\overline{D}$.
It is easy to see that this together with (1.3) yields
$$\displaystyle
\sup_{x\in\Omega\setminus\overline{D},\,\epsilon<\text{dist}\,(x,\partial\Omega)}
\Vert(w_1)_x\Vert_{L^2(\Omega\setminus\overline{D})}<\infty.
\tag {A.2}
$$
Note that we are considering general ${\cal G}$.
Using a similar argument for the proof of (19) in \cite{Iwave}, we have the estimate
$$\displaystyle
\Vert(w_2)_x\Vert_{L^2(\Omega\setminus\overline{D})}\le C\Vert G(\,,x)\Vert_{L^{\frac{4}{3}}(\partial D_d)}.
$$
Then assumption (1.3) and 
$$\displaystyle
\sup_{x\in\Bbb R^3}\Vert G(\,\cdot\,-x)\Vert_{L^{\frac{4}{3}}(\partial D_d)}<\infty,
$$
we obtain
$$\displaystyle
\sup_{x\in\Omega\setminus\overline{D},\,\epsilon<\text{dist}\,(x,\partial\Omega)}
\Vert(w_2)_x\Vert_{L^2(\Omega\setminus\overline{D})}<\infty.
\tag {A.3}
$$
From (A.2) and (A.3) we obtain the desired conclusion.

\noindent
$\Box$

\end{document}